\documentclass[12pt]{article}
\usepackage{amsmath,amsthm,amsfonts,amssymb,epsfig,verbatim,xcolor,constants,subcaption}
\usepackage[left=1in,top=1in,right=1in]{geometry}
\usepackage[normalem]{ulem}

\newconstantfamily{smc}{symbol=c}
\newconstantfamily{lgc}{symbol=C}

\newtheorem{thm}{Theorem}[section]
\newtheorem{lem}[thm]{Lemma}
\newtheorem{prop}[thm]{Proposition}

\newtheorem{df}[thm]{Definition}
\newtheorem{claim}[thm]{Claim}
\newtheorem{rem}{Remark}

\numberwithin{equation}{section}

\title{Estimates for the empirical distribution along a geodesic in first-passage percolation}
\author{Michael Damron \thanks{The research of M. D. was supported by an NSF CAREER grant and NSF grant DMS-2054559.} \\ \small{Georgia Tech}  \and Jack Hanson  \thanks{The research of J. H. was supported by NSF grant DMS-1954257.}\\ \small{City College, CUNY} \and Chris Janjigian \thanks{The research of C. J. is supported by NSF grant DMS-2125961.} \\ \small{Purdue University} \and Wai-Kit Lam\thanks{The research of W.-K. L. is supported by the National Science and Technology Council in Taiwan grant numbers 110-2115-M002-012-MY3 and  113-2115-M-002-009-MY3.} \\ \small{National Taiwan University} \and Xiao Shen \\ \small{University of Utah} }

\begin{document}
	
	\maketitle 
	\begin{abstract}
		In first-passage percolation, we assign i.i.d.~nonnegative weights $(t_e)$ to the nearest-neighbor edges of $\mathbb{Z}^d$ and study the induced pseudometric $T = T(x,y)$. In this paper, we focus on geodesics, or optimal paths for $T$, and estimate the empirical distribution of weights along them. We prove an upper bound for the expected number of edges with weight $\geq M$ in the union of all geodesics from $0$ to $x$ of the form $q(M) \mathbb{P}(t_e \geq M)|x|$, where $q(M) \leq e^{-cM}$. This shows that the tail of the expected empirical distribution along a geodesic is lighter than that of the original weight distribution by an exponential factor. We also give a lower bound for the expected minimal number of edges with weight $\geq M$ in any geodesic from $0$ to $x$ in terms of $\mathbb{P}(t_e \geq M)$ and $\mathbb{P}(t_e \in [M,2M])$. For example, these two imply that if $t_e$ has a power law tail of the form $\mathbb{P}(t_e \geq M) \sim M^{-\alpha}$, then the tail of the expected empirical distribution asymptotically lies between $e^{-CM \log M}$ and $e^{-cM}$. We also provide estimates for the expected number of edges in a geodesic with weight in a set $A$ for (a) arbitrary $A$, (b) $A$ an interval separated from the infimum of the support of $t_e$ and (c) $A=[0,a]$ for some $a \geq 0$. 
	\end{abstract}

\section{Introduction}

\subsection{Background and the model}
First-passage percolation (FPP) has been a important object of study in probability theory since it was originally introduced by Hammersley and Welsh \cite{HW65} nearly sixty years ago as a model of the flow of a fluid through a porous medium. It also arises naturally in a number of other settings, the most relevant of which for our purpose is in defining a random pseudo-metric on the vertices of a graph. For a recent account of research in FPP, we refer the reader to the survey \cite{ADH17}.

We will consider FPP defined as usual on the $d$-dimensional cubic lattice $\mathbb{Z}^d$ with the set of nearest-neighbor edges $\mathcal{E}^d$. Let $(t_e)_{e \in \mathcal{E}^d}$ be an i.i.d.~collection of nonnegative random variables with common distribution function $F$. A \underline{path} $\gamma$ between vertices $x$ and $y$ is a sequence $x=v_0,e_0,v_1,e_1, \dots, e_{n-1},v_n = y$ of alternating vertices and edges such that for each $i$, the edge $e_i$ has endpoints $v_i$ and $v_{i+1}$. The corresponding \underline{passage time of $\gamma$} is $T(\gamma) = \sum_{i=0}^{n-1} t_{e_i}$. The \underline{passage time between vertices $x$ and $y$} is
\[
T(x,y) = \inf_{\gamma : x \to y} T(\gamma),
\]
where the infimum is over all paths $\gamma$ from $x$ to $y$. 

A \underline{geodesic} from $x$ to $y$ is a path $\gamma$ from $x$ to $y$ such that $T(\gamma) = T(x,y)$. It is known \cite[Prop.~4.4]{ADH17} that almost surely, geodesics exist between all $x,y \in \mathbb{Z}^d$ under the assumption
\begin{equation}\label{eq: percolation_condition}
F(0) < p_c.
\end{equation}
Here $p_c$, the critical threshold for bond percolation in $d$ dimensions, is a number in $(0,1)$ that equals $1/2$ for $d=2$, and is decreasing in $d$. If $F$ has discontinuities, then almost surely, there exist vertices between which there are multiple geodesics. We write $\overline{\textsc{GEO}}(x,y)$ for the union
\[
\overline{\textsc{GEO}}(x,y) = \{e \in \mathcal{E}^d : e \text{ is in some geodesic from } x \text{ to } y\}.
\]
Also, we write $\underline{\textsc{GEO}}(x,y)$ for the intersection
\[
\underline{\textsc{GEO}}(x,y) = \{e \in \mathcal{E}^d : e \text{ is in every geodesic from } x \text{ to } y\}.
\]
We are interested in the empirical distribution of edge-weights on geodesics, so for any Borel set $A \subset \mathbb{R}$ and vertices $x,y$, we set
\begin{align*}
\overline{N}(x,y,A) &= \#\left\{e \in \overline{\textsc{GEO}}(x,y) : t_e \in A\right\} \text{ and}\\
\underline{N}(x,y,A) &= \#\left\{e \in \underline{\textsc{GEO}}(x,y) : t_e \in A\right\}.
\end{align*}
We also define
\begin{align*}
N_{\mathrm{max}}(x,y,A) &= \sup_\gamma \#\left\{e \in \gamma : t_e \in A\right\} \text{ and}\\
N_{\mathrm{min}}(x,y,A) &= \inf_\gamma \#\left\{e \in \gamma : t_e \in A\right\},
\end{align*}
where the supremum and infimum are over all geodesics $\gamma$ from $x$ to $y$. By definition, we have the inequalities
\[
\underline{N}(x,y,A) \leq N_{\mathrm{min}}(x,y,A) \leq N_{\mathrm{max}}(x,y,A) \leq \overline{N}(x,y,A).
\]

Our main goal is to give estimates for the expectations of the above four quantities. We seek bounds of the form $c(A)|x-y|$, where $c(A)$ is a constant depending only on $A$, and the Euclidean norm $|x-y|$ corresponds to the fact that geodesics typically have linear length. 
We now discuss previous work related to this goal. First, van den Berg and Kesten proved \cite[Rem.~2.15]{BK93} under the slighty stronger condition
\begin{equation}\label{eq: useful}
F(0) < p_c \text{ and } F(r) < \vec{p}_c,
\end{equation}
where $r = \inf\{x \in \mathbb{R} : F(x) > 0\}$ is the infimum of the support of the distribution of $t_e$, that if $A$ is any Borel set such that $\mathbb{P}(t_e \in A) > 0$, then $\mathbb{E}N_{\mathrm{min}}(0,n\mathbf{e}_1,A)/n$ is bounded away from 0 as $n \to \infty$. This statement can be formulated in terms of the empirical distribution of weights along a geodesic. If $\gamma$ is any geodesic connecting $0$ to $n\mathbf{e}_1$ with edges $e_1, e_2, \dots, e_k$ in order, then the empirical distribution of weights along $\gamma$ is defined as the distribution of a uniformly selected weight along $\gamma$:
\[
\mu_\gamma(A) = \frac{1}{k} \sum_{i=1}^k \mathbf{1}_{\{t_{e_i} \in A\}} \text{ for Borel } A \subset \mathbb{R}.
\]
The result of van den Berg and Kesten can be used, along with estimates on geodesic lengths, to show that under \eqref{eq: useful}, one has
\[
\mathbb{P}(t_e \in A) > 0 \Rightarrow \liminf_{n \to \infty} \mathbb{E}\min_\gamma \mu_\gamma (A) > 0,
\]
where the minimum is over all geodesics from $0$ to $n\mathbf{e}_1$. Furthermore, any limit point as $n \to \infty$ of the averaged empirical distributions $\mathbb{E}\mu_{\gamma_n}(\cdot )$ for $\gamma_n$ a (measurably chosen) geodesic from $0$ to $n \mathbf{e}_1$ is absolutely continuous relative to the distribution of $t_e$. Because of the relation between empirical distributions and edge counts, and for simplicity, we state all of our theorems in terms of the variables $\underline{N},\overline{N},N_{\min},$ and $N_{\max}$.

The second result on empirical distributions we mention, due to Damron-Hanson-Sosoe, is \cite[Thm.~6.6]{DHS15}: if \eqref{eq: percolation_condition} holds and $\mathbb{E}\min\{t_1, \dots, t_{2d}\}^\alpha < \infty$ for some $\alpha>0$ (the $t_i$ are i.i.d.~random variables distributed as $t_e$), then there exists $\Cl[lgc]{c: DHS}>0$ such that for all $x \in \mathbb{Z}^d$, we have
\begin{equation}\label{eq: DHS_sublinear}
\mathbb{E}N_{\mathrm{max}}(0,x,A) \leq \Cr{c: DHS} \mathbb{P}(t_e \in A)^{\frac{\alpha-1}{\alpha d}} |x|.
\end{equation}
(See Prop.~\ref{prop: DHS_improvement} below, where we improve this exponent to $1/d$ without a moment condition.) Next, Nakajima \cite{N22} studied the maximal weight of any edge on a geodesic, and determined its rate of growth depending on the tail of the distribution of $t_e$.  Bates \cite{B20} has shown, partially answering a question of Hoffman, that empirical distributions along geodesics from $0$ to $x$ have limits as $|x|$ grows, so long as the distribution of $t_e$ lies in various dense subsets of the set of all distributions. In a related exactly solvable model, vertex last-passage percolation with exponential weights, Martin-Sly-Zhang have shown \cite{MSZ21} convergence of empirical distributions, and even convergence of the environment seen from the geodesic.

The current study was motivated by the question of whether the exponent $(\alpha-1)/\alpha d$ (or $1/d$) from \eqref{eq: DHS_sublinear} (or Prop.~\ref{prop: DHS_improvement} below) is optimal. In Thm.~\ref{thm: upper_bound_upper_tail}, we consider sets of the form $A = (M,\infty)$ and prove that, in this case, the right side of \eqref{eq: DHS_sublinear} can be improved to $e^{-cM} \mathbb{P}(t_e \in A)|x|$. Furthermore, the left side can be replaced by $\mathbb{E}\overline{N}(0,x,A)$. That is, not only can the exponent $(\alpha-1)/\alpha d$ be improved to 1 but, depending on the tail of the original edge-weight distribution, the exponential factor $e^{-cM}$ can make the bound much smaller. For example, if $t_e$ has a power-law tail, then the right side is exponential in $M$. Further, for general $t_e$, the bound shows that the empirical distribution along a geodesic has a lighter tail (sometimes significantly lighter) than that of the distribution of $t_e$. In Prop.~\ref{prop: separation_prop}, we consider sets $A$ which are bounded away from the infimum of the distribution of $t_e$. Specifically, if we set $r = \inf\{x : F(x) > 0\}$ and fix any $b > r$, then for any interval $A = (c,d] \subset [b,\infty)$, the right side of \eqref{eq: DHS_sublinear} can be replaced by $C\mathbb{P}(t_e \in A) |x|$, where $C$ depends on $b$. In other words, the exponent for such $A$ can also be improved to 1.

In this paper, we also give some lower bounds for the quantities $N_{\textrm{min}}(0,x,A)$ and $\underline{N}(0,x,A)$. The main such result is Thm.~\ref{thm: BK_argument}, which estimates $\mathbb{E}N_{\textrm{min}}(0,x,[M,\infty))$ for $M$ sufficiently large and provides a lower bound in terms of $\mathbb{P}(t_e \in [M,2M])$ and $\mathbb{P}(t_e \geq M)$. For example, in the case that $t_e$ has a power-law tail, the bound is of the form $e^{-M\log M}|x|$, and this is close to the exponential upper bound mentioned above; see Rem.~\ref{rem: BK_examples}. A simpler inequality is also given for $A$ of the form $[0,a]$ in Prop.~\ref{prop: FKG_lower_bound} of the form $\mathbb{E}\underline{N}(0,x,A) \geq c\mathbb{P}(t_e \in A)|x|$.

\subsection{Main results}

\subsubsection{Bounds for the right tail} 
Our first main result gives an upper bound for $\overline{N}(0,x,A)$ in the case that $A = (M,\infty)$ for some $M \geq 0$. Intuitively it states that the expected empirical distribution of edge-weights along a geodesic has a tail which is lighter than that of the original edge-weight distribution by a factor of $q(M,x)$, a quantity that has exponential decay in the sense that $\limsup_{|x|\to\infty}q(M,x)/|x| \leq e^{-cM}$. In other words, for large $M$, we have 
\[
\limsup_{|x| \to \infty} \frac{\mathbb{E}\overline{N}(0,x,(M,\infty))}{|x|} \leq C (1-F(M))e^{-cM}.
\]
In the statement of Thm.~\ref{thm: upper_bound_upper_tail} below, we use the convention that $(1-F(M))/F(M) = \infty$ when $F(M)=0$.

\begin{thm}\label{thm: upper_bound_upper_tail}
Suppose that \eqref{eq: percolation_condition} holds. There exist $\Cl[smc]{c: L2}, \Cl[lgc]{c: L1}>0$ such that for all $x \in \mathbb{Z}^d$ and all $M\geq 0$,
\[
\mathbb{E}\overline{N}(0,x,(M,\infty)) \leq \frac{1-F(M)}{F(M)} q(M,x),
\]
where $q(M,x)$ is the minimum of the following two quantities:
\begin{enumerate}
\item $\Cr{c: L1}(1+ (1-F(M/4))^{2-2/d}|x|)$ and
\item $\Cr{c: L1}(1+M^d)(1+e^{-\Cr{c: L2} M}|x|)$.
\end{enumerate}
\end{thm}
The proof of Thm.~\ref{thm: upper_bound_upper_tail} uses a version of the observation from \cite{N22} that if $e$ is in a geodesic from $0$ to $x$ but $t_e \geq M$, then we cannot reroute the geodesic around $e$ to take a shortcut that has total passage time less than $M$. Instead of applying this idea to the event that $e$ is in a geodesic, we apply it in Claim~\ref{claim: good_claim} to the event $\{D_e > M\}$, where $D_e$ is a variable that measures the stability of the geodesic relative to perturbations in the weight $t_e$. To show that typical edges are likely to have many shortcuts around them (these are the ``good'' edges of Def.~\ref{def: good}), in Lem~\ref{lem: good_lemma} we use techniques developed for studying the chemical distance in supercritical bond percolation. Last, we estimate the maximal number of edges without shortcuts that exist along any path from $0$ to $x$ using an inequality proved in Appendix~\ref{sec: appendix} (in Prop.~\ref{prop: lattice_animals}) for the greedy lattice animal problem under a finite-dependence condition.

In the next remark, and in many of the proofs, we use the notation
\begin{equation}\label{eq: geo_long_def}
\textsc{GEO}(0,x) = \text{the vertex self-avoiding geodesic from }0 \text{ to }x\text{ with the most edges},
\end{equation}
and $\#\textsc{GEO}(0,x)$ is the number of edges in $\textsc{GEO}(0,x)$. If there are multiple such geodesics with the most edges, we choose the first one in some deterministic ordering of all finite paths.

\begin{rem}
The expression on the right side in Thm.~\ref{thm: upper_bound_upper_tail} is infinite if $M$ is less than the infimum of the distribution of $t_e$. Indeed, for such $M$ we should not expect a bound that is linear in $x$ because the quantity on the left is $\mathbb{E}\# \overline{\textsc{GEO}}(0,x)$, which can be of order $|x|^d$ when the $t_e$ are constant random variables. Instead, if we consider $N_{\mathrm{max}}(0,x,[M,\infty))$, then for $M$ such that $F(M-) \geq 1/2$, Thm.~\ref{thm: upper_bound_upper_tail} implies
\[
\mathbb{E}N_{\mathrm{max}}(0,x,[M,\infty)) \leq 2 (1-F(M-)) q(M,x),
\]
whereas for $M$ such that $F(M-) \leq 1/2$, \eqref{eq: precise_geodesic_moment_bound} gives for some $\Cl[lgc]{c: tempman_numero_uno}>0$ and all nonzero $x$
\[
\mathbb{E}N_{\mathrm{max}}(0,x,[M,\infty)) \leq \mathbb{E}\#\textsc{GEO}(0,x) \leq \Cr{c: tempman_numero_uno}|x|.
\]
Putting these together shows that for some $\Cl[lgc]{c: L6}>0$ and all $M \geq 0$, 
\[
\mathbb{E}N_{\mathrm{max}}(0,x,[M,\infty)) \leq \Cr{c: L6} (1-F(M-)) q(M,x).
\]
\end{rem}

\begin{rem}
For heavy-tailed distributions, the second term in the minimum defining $q(M,x)$ is smaller than the first as $M \to \infty$. For example, if $1-F(M)$ is of order $M^{-\alpha}$, then the first term is of order $M^{-\alpha(2-2/d)} |x|$, whereas the second is of order $e^{-cM}|x|$. For light-tailed distributions, however, the first term can be smaller. Taking $1-F(M)$ instead to be of order $e^{-M^\alpha}$, the second term is still of order $e^{-cM}|x|$, while the first is of order $e^{-cM^\alpha}|x|$.
\end{rem}

In the second main result, we give lower bounds for the number of edges with weight at least $M$ in geodesics. To do this, we will need to assume van den Berg and Kesten's condition \eqref{eq: useful}.
In the absence of \eqref{eq: useful}, long geodesics can contain all but the first few edges with weight equal to $r$, and so inequalities like those in Thm.~\ref{thm: BK_argument} cannot hold.

\begin{thm}\label{thm: BK_argument}
Suppose that \eqref{eq: useful} holds. There exist constants $\Cl[smc]{c: BK_lower_constant}, \Cl[lgc]{c: BK_upper_constant} >0$ such that for all $x \in \mathbb{Z}^d$ and all $M$ sufficiently large, we have
\[
\mathbb{E} N_{\mathrm{min}}(0,x,[M,\infty)) \geq \mathbb{P}(t_e \in [M,2M])^{2d} \mathbb{P}(t_e \geq M)^{\Cr{c: BK_upper_constant}M} \left( \frac{\Cr{c: BK_lower_constant}}{M} |x|- \Cr{c: BK_upper_constant}\right).
\]
\end{thm}

The proof of Thm.~\ref{thm: BK_argument} uses tools from the van den Berg-Kesten argument \cite{BK93} (see also \cite{N22}, where these were used to give a lower bound for the maximal edge-weight in a geodesic). Roughly speaking, we show that any geodesic from $0$ to $x$ must cross many annuli with inradius $cM$ and outradius $CM$ in such a way that the crossing does not have a passage time that is too small. These annuli correspond to the ``black'' cubes of Def.~\ref{def: black_cube}. In Prop.~\ref{prop: modification}, we do a edge-modification (``resampling'') argument to replace the geodesic from its first entrance of a portion of the annulus to its last exit with a deterministic ``mostly oriented'' detour with low passage time. These detours are constructed in such a way that in the resulting resampled environment, there can be only one fast detour, and all edges on the boundary of the detour have weight above $M$. Futhermore, the geodesic in the resampled environment must intersect the detour, so it must take an edge with weight above $M$.

The lower bound of Thm.~\ref{thm: BK_argument} implies an estimate of the form
\begin{equation}\label{eq: BK_asymptotic_form}
\liminf_{|x| \to \infty} \frac{\mathbb{E}N_{\mathrm{min}}(0,x,[M,\infty))}{|x|} \geq \mathbb{P}(t_e \in [M,2M])^{2d} \mathbb{P}(t_e \geq M)^{\Cr{c: BK_upper_constant}M}.
\end{equation}
Elementary modifications of the proof of Thm.~\ref{thm: BK_argument} can provide analogous lower bounds for more general sets $A$ than just  $A = [M,\infty)$. As one example, if we change items 3 and 4 in the definition of $A_{u,v}^\ast$ (above \eqref{eq: A_bound}) to ``$t_e^\ast \in [M,2M)$ for $e \in E_3(u,v) \cup E_4(u,v)$'' we produce 
\[
\mathbb{E}N_{\mathrm{min}}(0,x,[M,2M)) \geq \mathbb{P}(t_e \in [M,2M))^{\Cr{c: BK_upper_constant}M} \left( \frac{\Cr{c: BK_lower_constant}}{M} |x| - \Cr{c: BK_upper_constant}\right),
\]
and therefore the following variant of \eqref{eq: BK_asymptotic_form}:
\[
\liminf_{|x| \to \infty} \frac{\mathbb{E}N_{\mathrm{min}}(0,x,[M,2M))}{|x|} \geq \mathbb{P}(t_e \in [M,2M))^{\Cr{c: BK_upper_constant}M}.
\]
If we sum this inequality over $M,2M,4M, \dots$ and use a simple limiting argument, we obtain
\[
\liminf_{|x| \to \infty} \frac{\mathbb{E}N_{\mathrm{min}}(0,x,[M,\infty))}{|x|} \geq \sum_{k=0}^\infty \mathbb{P}(t_e \in [2^{k-1}M,2^kM))^{2^{k-1}\Cr{c: BK_upper_constant}M}.
\]
For some distributions of $t_e$, like those which have gaps that imply $\mathbb{P}([M,2M]) = 0$, this improves on the estimate of Thm.~\ref{thm: BK_argument}.

\begin{rem}\label{rem: BK_examples}
If $t_e$ has a power law tail of the form $1-F(M) \sim M^{-\alpha}$, then the lower bound in Thm.~\ref{thm: BK_argument} is of the form $e^{-CM \log M}|x|$, which is close to the upper bound $e^{-cM}|x|$ that comes from Thm.~\ref{thm: upper_bound_upper_tail}. In the case of a lighter tail of the form $1-F(M) \sim e^{-M^\alpha}$, the lower bound of Thm.~\ref{thm: BK_argument} is of order $e^{-CM^{\alpha+1}}|x|$, compared to the upper bound from Thm.~\ref{thm: upper_bound_upper_tail} of the form $e^{-c\max\{M,M^{\alpha}\}}|x|$. 
\end{rem}

\subsection{Bounds for general sets $A$}

In this section, we give bounds for the number of edges in geodesics whose weights lie in some Borel set $A$, not necessarily of the form $(M,\infty)$. The first is an improvement on \cite[Thm.~6.6]{DHS15}, listed in \eqref{eq: DHS_sublinear}, on the number of edges in a geodesic with weight in an arbitrary Borel set $A$. Specifically, it removes the moment condition $\mathbb{E}\min\{t_1, \dots, t_{2d}\}^\alpha < \infty$ (where the $t_i$ are i.i.d.~edge-weights and $\alpha>0$), and improves the exponent $(\alpha-1)/\alpha d$ to $1/d$. 

\begin{prop}\label{prop: DHS_improvement}
Suppose that \eqref{eq: percolation_condition} holds. There exists $\Cl[lgc]{c: DHS_improvement}>0$ such that for all $x \in \mathbb{Z}^d$ and all Borel sets $A \subset \mathbb{R}$, we have
\[
\mathbb{E}N_{\mathrm{max}}(0,x,A) \leq \Cr{c: DHS_improvement} \mathbb{P}(t_e \in A)^{\frac{1}{d}}|x|.
\]
\end{prop}

The second result is an upper bound for sets $A$ that have some positive distance from the infimum $r$ of the distribution of $t_e$. In this case, it improves the exponent $1/d$ from Prop.~\ref{prop: DHS_improvement} to $1$. It also coincides with our intuition that the empirical distribution along a geodesic should be lighter than that of the original edge-weight distribution, with more mass possibly concentrated near $r$. For such $A$, the result implies the estimate
\[
\limsup_{|x| \to \infty} \frac{\mathbb{E}\overline{N}(0,x,A)}{|x|} \leq \Cr{c: general_trick_bound} \mathbb{P}(t_e \in A).
\]
Prop.~\ref{prop: separation_prop} is stated for intervals $A = (c,d]$ with $b \leq c \leq d$, but because the constant prefactor $\Cr{c: general_trick_bound}$ is uniform over such choices of $c,d$ it extends by a standard monotone class argument to all Borel subsets $A$ of $(b,\infty)$.
\begin{prop}\label{prop: separation_prop}
Suppose that \eqref{eq: percolation_condition} holds and let $b$ be such that $F(b) > F(r)$. There exists $\Cl[lgc]{c: general_trick_bound} >0$ such that for all $c,d$ with $b \leq c \leq d$, we have
\[
\mathbb{E}\overline{N}(0,x,(c,d]) \leq \Cr{c: general_trick_bound} \mathbb{P}(t_e \in (c,d]) |x| \text{ for all } x \in \mathbb{Z}^d.
\]
\end{prop}

For $A$ of the form $[0,a]$, (with $a \geq 0$), we have the following estimate.

\begin{prop}\label{prop: FKG_lower_bound}
Suppose that \eqref{eq: useful} holds. There exists $\Cr{c: easy_lower_bound}$ such that for all $x \in \mathbb{Z}^d$ and all $a \geq 0$, we have
\[
\mathbb{E}\underline{N}(0,x,[0,a]) \geq \Cl[smc]{c: easy_lower_bound} \mathbb{P}(t_e \in [0,a]) |x|.
\]
\end{prop}

We do not know if the lower bound $c \mathbb{P}(t_e \in A)|x|$ from the last proposition is the largest possible for general distributions. In the following Bernoulli example with $A = \{0\} = [0,a]$ for $a = 0$, we can find a larger lower bound $c\mathbb{P}(t_e \in A)^{1/d}|x|$.

\medskip
\noindent
{\bf Example.} (Bernoulli FPP) 
For $p \in (0,1)$, take $(t_e)$ to be i.i.d.~with 
\[
\mathbb{P}(t_e = 0) = p = 1-\mathbb{P}(t_e=1)
\]
and $x_n = n \sum_{i=1}^d \mathbf{e}_i = (n, \dots, n)$. 
If $\gamma$ is any geodesic from $0$ to $x_n$, then
\[
\#\{e \in \gamma : t_e = 0\} = \# \text{ edges in } \gamma - \# \{e \in \gamma : t_e = 1\} \geq nd - T(\gamma) = nd-T(0,x_n),
\]
and so
\begin{equation}\label{eq: oriented_comparison}
\mathbb{E}N_{\mathrm{min}} (0,x_n, \{0\}) \geq \mathbb{E} \left( nd - T(0,x_n) \right).
\end{equation}

To estimate $T(0,x_n)$, we construct a (probably non-optimal) path. Let $n_p = \lceil p^{-1/d}\rceil$ and define $E_p$ to be the event that there is at least one edge $e$ contained in the box $[0,n_p]^d$ such that $t_e=0$. Then there exists $\Cl[smc]{c: poisson}>0$ such that for all $p \in (0,1)$, we have $\mathbb{P}(E_p) \geq \Cr{c: poisson}$. On the event $E_p$, let $e^\ast$ be the first edge $e$ in some deterministic ordering in this box with $t_e=0$ and define an oriented (moving only in the positive coordinate directions) path $\gamma^\ast$ from $0$ to $x_{n_p}$ by starting at 0, moving to the near endpoint of $e^\ast$, crossing $e^\ast$, and then moving to $x_{n_p}$. By construction, $T(\gamma^\ast) \leq n_p d-1$, hence
\begin{equation}\label{eq: poisson_lower_bound}
\mathbb{E}(n_pd - T (0,x_{n_p})) \geq \mathbb{E}(n_p d - T (\gamma^\ast))\mathbf{1}_{E_p} \geq \Cr{c: poisson}.
\end{equation}
For any $n \geq n_p$, we can write $n = k n_p + r$ where $k = \lfloor n / n_p \rfloor$ and estimate
\[
\mathbb{E}T (0,x_n) \leq \sum_{j=0}^{k-1} \mathbb{E}T (x_{jn_p}, x_{(j+1)n_p}) + rd = k\mathbb{E}T (0,x_{n_p}) + rd,
\]
and so, using \eqref{eq: oriented_comparison} and \eqref{eq: poisson_lower_bound}, we have
\[
\mathbb{E}N_{\mathrm{min}}(0,x_n,\{0\}) \geq nd - k \mathbb{E}T(0,x_{n_p}) - rd = k (n_p d - \mathbb{E}T (0,x_{n_p})) \geq \Cr{c: poisson} k,
\]
which is of order $np^{1/d}$.

\subsection{Outline of the paper}
The rest of the paper consists of proofs of the main results. We begin in Sec.~\ref{sec: tail_bound} with a tail bound on the length of geodesics under only assumption~\eqref{eq: percolation_condition}. This result will be used in some of the subsequent proofs. Next, Thm.~\ref{thm: upper_bound_upper_tail} is proved in Sec.~\ref{sec: upper_bound_upper_tail}, and Thm.~\ref{thm: BK_argument} is proved in Sec.~\ref{sec: BK_argument}. Sec.~\ref{sec: general_A} contains proofs of Prop.~\ref{prop: DHS_improvement}, Prop.~\ref{prop: separation_prop}, and Prop.~\ref{prop: FKG_lower_bound}. In Appendix~\ref{sec: appendix}, we prove an upper bound for the finitely-dependent greedy lattice animal problem which extends \cite[Lem.~6.8]{DHS15} and is important in the proof of Thm.~\ref{thm: upper_bound_upper_tail}.

\section{Tail bound for geodesic length}\label{sec: tail_bound}

Recall the definition of $\textsc{GEO}(0,x)$ from \eqref{eq: geo_long_def}. In this section we will extend the argument from \cite[Thm.~4.9]{ADH17} to prove the following estimate. 

\begin{prop}\label{prop: improved_tail_bound}
Suppose that \eqref{eq: percolation_condition} holds. There exist $\Cl[smc]{c: improved_small_constant}, \Cl[lgc]{c: improved_big_constant}>0$ such that for $\lambda \geq \Cr{c: improved_big_constant}$ and all nonzero $x \in \mathbb{Z}^d$, we have
\begin{equation}\label{eq: my_new_tail_bound}
\mathbb{P}\left( \#\textsc{GEO}(0,x) \geq \lambda |x|\right) \leq  \exp\left( - \Cr{c: improved_small_constant} \left( \lambda |x|\right)^{\frac{1}{d}}\right).
\end{equation}
In particular, for $p \geq 1$ and all nonzero $x$, we have
\begin{equation}\label{eq: precise_geodesic_moment_bound}
\mathbb{E}\left( \#\textsc{GEO}(0,x) \right)^p \leq \left( \Cr{c: improved_big_constant} |x|\right)^p + \frac{\Gamma(dp+1)}{\Cr{c: improved_small_constant}^{dp}}.
\end{equation}
\end{prop}

To prove this, we consider an embedded percolation model. Let $M_0>0$ be such that $F(M_0) > p_c$. We call an edge \underline{open} if $t_e \leq M_0$ and \underline{closed} otherwise. A path is said to be open if all its edges are open. An \underline{open cluster} is a maximal set of vertices any two of which are connected by an open path. By choice of $M_0$, almost surely there exists a (unique) infinite open cluster, which we denote by $I$. Write $B(n) = [-n,n]^d$ and its internal vertex boundary $\partial_i B(n) = \{x \in B(n) : \exists y \in B(n)^c \text{ such that } |x-y| = 1\}$. To ensure that geodesics intersect $I$, we use \cite[Lem.~6.3]{ADH14}:

\begin{lem}\label{lem: gamma_lemma}
Let $A_n$ be the event that every path from $0$ to $\partial_i B(n)$ intersects $I$. There exists $p_0 \in (p_c,1)$ such that if $F(M_0) \geq p_0$, then for some $\Cl[smc]{c: book_chemical_constant}>0$,
\[
\mathbb{P}(A_n) \geq 1 - e^{-\Cr{c: book_chemical_constant}n} \text{ for all } n.
\]
\end{lem}
\noindent
For the rest of this section, we fix $M_0$ such that $F(M_0) \geq p_0$.

We will also need a result from \cite[Rem.~1.7]{DN22} which follows from \cite[Eq.~(4.49)]{AP96}. Here we use the notation $d_I(x,y)$ to mean the smallest number of edges in any open path from $x$ to $y$. If there is no such path, we set $d_I(x,y) = \infty$.
\begin{lem}\label{lem: AP_lem}
There exist $\Cl[smc]{c: AP_small_constant}, \Cl[lgc]{c: AP_large_constant}>0$ such that for all $x \in \mathbb{Z}^d$ and all $\ell \geq \Cr{c: AP_large_constant}|x|$, we have 
\[
\mathbb{P}(\ell \leq d_I(0,x) < \infty) \leq e^{-\Cr{c: AP_small_constant}\ell}.
\]
\end{lem}

\begin{proof}[Proof of Prop.~\ref{prop: improved_tail_bound}]
First, \eqref{eq: precise_geodesic_moment_bound} follows from \eqref{eq: my_new_tail_bound} by integrating. Namely, we write
\begin{align*}
\mathbb{E}\left( \#\textsc{GEO}(0,x) \right)^p & = |x|^p \int_0^\infty p \lambda^{p-1}\mathbb{P}(\#\textsc{GEO}(0,x) \geq \lambda |x|)~\text{d}\lambda \\
&\leq (\Cr{c: improved_big_constant} |x|)^p + |x|^p \int_0^\infty p\lambda^{p-1}\exp\left( - \Cr{c: improved_small_constant} \left( \lambda |x|\right)^{\frac{1}{d}}\right)~\text{d}\lambda \\  
&= (\Cr{c: improved_big_constant} |x|)^p + \frac{dp}{\Cr{c: improved_small_constant}^{dp}} \int_0^\infty u^{dp-1}e^{-u}~\text{d}u.
\end{align*}
This is the claimed estimate in \eqref{eq: precise_geodesic_moment_bound}.

To begin the proof of \eqref{eq: my_new_tail_bound}, we let $x \in \mathbb{Z}^d$ be nonzero and, for $\lambda \geq 2^{d+2}$, let $B_1$ be the cube of sidelength $(\lambda |x|/2^{d+2})^{1/d}$ centered at the origin, with the cube $B_2$ defined as the translation $B_2 = B_1 + x$. Consider an outcome for which $\#\textsc{GEO}(0,x) \geq \lambda |x|$. Because 
\begin{equation}\label{eq: B_i_volume_count}
\# \text{ vertices in } (B_1 \cup B_2) \leq 2 \left( 1 + \left( \frac{\lambda |x|}{2^{d+2}}\right)^{\frac{1}{d}}\right)^d \leq \frac{\lambda |x|}{2},
\end{equation}
$\textsc{GEO}(0,x)$ must contain a vertex, say $u$, that is not in $B_1 \cup B_2$. Let $u_1$ be the first vertex of $\textsc{GEO}(0,x)$ in $B_1^c$ and let $u_2$ be the last vertex of $\textsc{GEO}(0,x)$ in $B_2^c$. As we follow $\textsc{GEO}(0,x)$ from $0$ to $x$, the vertex $u_1$ appears first, then $u$, then $u_2$, although some of these may be equal. Let $\pi_1$ be the subpath of $\textsc{GEO}(0,x)$ from $0$ to $u_1$, and let $\pi_2$ be the path obtained by following $\textsc{GEO}(0,x)$ in reverse from $x$ to $u_2$. Then $\pi_1$ and $\pi_2$ are edge-disjoint.

Let $A(1)$ be the event $A_{\lceil (\lambda |x|/2^{d+2})^{1/d}/2\rceil}$ (this is the event $A_n$, where $n$ is chosen so that $u_1 \in \partial_i B(n)$) and let $A(2)$ be the event $A(1)$ naturally translated from $0$ to $x$. That is, the configurations $(t_e)$ in $A(2)$ are exactly those for which $(t_{e-x})$ is in $A(1)$. By Lem.~\ref{lem: gamma_lemma}, we have
\begin{align}
\mathbb{P}\left( \#\textsc{GEO}(0,x) \geq \lambda |x|\right) &\leq 2\exp\left( -\frac{\Cr{c: book_chemical_constant}}{2}\left( \frac{\lambda |x|}{2^{d+2}}\right)^{\frac{1}{d}}  \right) \nonumber \\
&+ \mathbb{P}\left(A(1) \cap A(2),~\#\textsc{GEO}(0,x) \geq \lambda |x|\right). \label{eq: pasta_macaroni_with_sauce}
\end{align}
For any outcome in the event appearing in \eqref{eq: pasta_macaroni_with_sauce}, the path $\pi_1$ must intersect a vertex $v_1 \in I$ and $\pi_2$ must intersect a vertex $v_2 \in I$. Because $\textsc{GEO}(0,x)$ is a geodesic, its subpath $\Gamma$ from $v_1$ to $v_2$ must also be a geodesic. By \eqref{eq: B_i_volume_count}, the number of vertices in $\Gamma$ is
\[
\# \Gamma \geq \#\textsc{GEO}(0,x) - \frac{\lambda |x|}{2} \geq \frac{\lambda |x|}{2}.
\]
We aim to show that the event appearing in \eqref{eq: pasta_macaroni_with_sauce} is unlikely, and we do this by showing that (a) because $\Gamma$ has many vertices ($\lambda$ is large), its passage time $T(v_1,v_2)$ should be large and (b) because $v_1$ and $v_2$ are in the infinite cluster, $T(v_1,v_2)$ should not be too large.

For (a), we need Kesten's lemma from \cite[Prop.~5.8]{aspects}. 
\begin{lem}\label{lem: kesten_lemma}
There exist $a,\Cl[smc]{c: kesten_lemma_constant}>0$ such that for every real $n \geq 1$, we have
\[
\mathbb{P}\left( \exists \text{ self-avoiding path } \sigma \text{ from } 0 \text{ with } \#\sigma \geq n \text{ but } T(\sigma) < an\right) \leq e^{-\Cr{c: kesten_lemma_constant}n}.
\]
\end{lem}
To use this in \eqref{eq: pasta_macaroni_with_sauce}, let $A(3)$ be the event that any self-avoiding path $\sigma$ starting at a vertex in the box $B_1$ with at least $\lambda |x|/2$ many vertices has $T(\sigma) \geq a \lambda |x|/2$.
By a union bound, \eqref{eq: B_i_volume_count}, and Lem.~\ref{lem: kesten_lemma}, we have
\[
\mathbb{P}(A(3)^c) \leq \frac{\lambda |x|}{2} \exp\left(- \Cr{c: kesten_lemma_constant} \frac{\lambda |x|}{2} \right),
\]
and so we can combine this with \eqref{eq: pasta_macaroni_with_sauce} to obtain
\begin{align}
\mathbb{P}(\#\textsc{GEO}(0,x) \geq \lambda |x|) &\leq 2 \exp\left( - \frac{\Cr{c: book_chemical_constant}}{2} \left( \frac{\lambda |x|}{2^{d+2}}\right)^{\frac{1}{d}} \right) + \frac{\lambda |x|}{2} \exp\left( - \Cr{c: kesten_lemma_constant}\frac{\lambda |x|}{2}\right) \nonumber \\
&+ \mathbb{P}\left(A(1) \cap A(2) \cap A(3),~ \#\textsc{GEO}(0,x) \geq \lambda |x|\right). \label{eq: hamburger_helper_texmex}
\end{align}
For any outcome in the event appearing in \eqref{eq: hamburger_helper_texmex}, we may find vertices $v_1,v_2 \in I$ as detailed above, and the geodesic $\Gamma$ connecting them. Because of occurrence of $A(3)$, we must have
\[
T(v_1,v_2) = T(\Gamma) \geq a \frac{\lambda |x|}{2}.
\] 
However, as $v_1,v_2$ are in the infinite cluster, the passage time between them should not be too large (item (b) mentioned above): there is a path connecting with $d_I(v_1,v_2)$ many edges, all of which have weight $\leq M_0$. Therefore
\[
M_0 d_I(v_1,v_2) \geq T(v_1,v_2) \geq a \frac{\lambda |x|}{2}.
\]
We have now shown that
\begin{equation}\label{eq: final_prob_bound}
\eqref{eq: hamburger_helper_texmex} \leq \mathbb{P}\left( \exists~v_1 \in B_1,~v_2 \in B_2 \text{ such that } \frac{a\lambda |x|}{2M_0} \leq d_I(v_1,v_2) < \infty \right).
\end{equation}
We will use Lem.~\ref{lem: AP_lem} to estimate this probability.

To apply Lem~\ref{lem: AP_lem}, we would like to have 
\begin{equation}\label{eq: we_want_to_apply_AP}
a\lambda |x|/(2M_0) \geq \Cr{c: AP_large_constant} |v_1-v_2|.
\end{equation}
To guarantee this, we increase $\lambda$ a bit more (so far it was required to be $\geq 2^{d+2}$ at the beginning of the proof) so that it further satisfies
\begin{equation}\label{eq: more_lambda_conditions}
\lambda \geq \frac{4\Cr{c: AP_large_constant}M_0}{a} \text{ and } \lambda \geq \left( \frac{8\Cr{c: AP_large_constant}M_0 d}{a} \right)^{\frac{d}{d-1}}.
\end{equation}
For $v_1 \in B_1$ and $v_2 \in B_2$, we have
\[
|v_1-v_2| \leq |x| + |v_1| + |v_2 - x| \leq |x| + 2 \max_{v \in B_1} |v| 
\leq |x| + 2d \left( \frac{\lambda |x|}{2^{d+2}}\right)^{\frac{1}{d}}.
\]
The first inequality of \eqref{eq: more_lambda_conditions} ensures that $\Cr{c: AP_large_constant} |x| \leq a\lambda |x|/(4M_0)$ and the second ensures that $2d\Cr{c: AP_large_constant} (\lambda |x| / 2^{d+2})^{1/d} \leq a\lambda |x|/(4M_0)$. In total, our conditions \eqref{eq: more_lambda_conditions} on $\lambda$ imply that for $v_1 \in B_1$, $v_2 \in B_2$, we have \eqref{eq: we_want_to_apply_AP}. Therefore we can return to \eqref{eq: final_prob_bound} and apply a union bound along with Lem.~\ref{lem: AP_lem} and \eqref{eq: B_i_volume_count} to obtain
\begin{align*}
\eqref{eq: hamburger_helper_texmex} \leq \sum_{v_1 \in B_1, v_2 \in B_2} \mathbb{P}\left( \frac{a \lambda |x|}{2M_0} \leq d_I(v_1,v_2) < \infty\right) &\leq (\# \text{vertices in } B_1)^2 \exp\left( - \Cr{c: AP_small_constant} \frac{a \lambda |x|}{2M_0} \right) \\
&\leq \left( \frac{\lambda |x|}{2}\right)^2 \exp\left( - \Cr{c: AP_small_constant} \frac{a \lambda |x|}{2M_0} \right).
\end{align*}
If we place this back in the equation that contains \eqref{eq: hamburger_helper_texmex}, we produce
\begin{align*}
\mathbb{P}(\#\textsc{GEO}(0,x) \geq \lambda |x|) &\leq 2 \exp\left( - \frac{\Cr{c: book_chemical_constant}}{2} \left( \frac{\lambda |x|}{2^{d+2}}\right)^{\frac{1}{d}} \right) + \frac{\lambda |x|}{2} \exp\left( - \Cr{c: kesten_lemma_constant}\frac{\lambda |x|}{2}\right)  \\
&+ \left( \frac{\lambda |x|}{2}\right)^2 \exp\left( - \Cr{c: AP_small_constant} \frac{a \lambda |x|}{2M_0} \right). 
\end{align*}
This implies the statement of the proposition.
\end{proof}

\section{Proof of Thm.~\ref{thm: upper_bound_upper_tail}}\label{sec: upper_bound_upper_tail}

To prove Thm.~\ref{thm: upper_bound_upper_tail}, we may suppose that $F(M) \in (0,1)$. To bound $\mathbb{E} \overline{N}(0,x,(M,\infty))$, we write
\begin{equation}\label{eq: initial_inequality}
\mathbb{E}\overline{N}(0,x,(M,\infty)) = \sum_e \mathbb{P}\left(e \in \overline{\textsc{GEO}}(0,x), t_e > M\right).
\end{equation}
For a fixed edge $e$, write our edge-weight configuration as a pair $(t_e,t_e^c)$, where $t_e^c = (t_f)_{f \neq e}$ is the collection of weights on edges not equal to $e$. Next, define
\begin{equation}\label{eq: D_e_def}
D_e = \sup\left\{s \geq 0 : e \in \overline{\textsc{GEO}}(0,x) \text{ in } (s,t_e^c)\right\}
\end{equation}
for the largest value of $s$ for which $e$ lies in a geodesic from $0$ to $x$ in the configuration in which we replace the value of $t_e$ by $s$. If the set above is empty, we set $D_e=-\infty$. We observe that $D_e$ is a function of $t_e^c$ (not $t_e$), and therefore is independent of $t_e$. Furthermore, by monotonicity of the event $\{e \in \overline{\textsc{GEO}}(0,x) \text{ in } (s,t_e^c)\}$ in $s$, we obtain
\begin{equation}\label{eq: overline_equivalence}
t_e \leq D_e \Leftrightarrow e \in \overline{\textsc{GEO}}(0,x).
\end{equation}
Last, we have
\[
t_e < D_e \Rightarrow e \in \underline{\textsc{GEO}}(0,x).
\]
The reason is that if $t_e < D_e$ but there is a geodesic $\Gamma_1$ from $0$ to $x$ that does not contain $e$, then in the configuration $\left( \frac{t_e+D_e}{2},t_e^c\right)$, $\Gamma_1$ would have strictly smaller passage time than any path from $0$ to $x$ that contains $e$, so $e$ would not be in $\overline{\textsc{GEO}}(0,x)$, a contradiction.

With this notation, we estimate \eqref{eq: initial_inequality} using independence and \eqref{eq: overline_equivalence} as
\begin{align}
&\mathbb{E}\overline{N}(0,x,(M,\infty)) \nonumber\\
=~& \sum_e \mathbb{P}(t_e \leq D_e, t_e > M)\nonumber \\
\leq~& \sum_e \mathbb{P}(D_e > M, t_e > M) \nonumber \\
=~& (1-F(M)) \sum_e \mathbb{P}(D_e > M) \nonumber\\
=~& \frac{1-F(M)}{F(M)} \sum_e \mathbb{P}(t_e \leq M, D_e > M) \nonumber\\
\leq~& \frac{1-F(M)}{F(M)} \sum_e \mathbb{P}(e \in \underline{\textsc{GEO}}(0,x), D_e > M) \nonumber\\
=~& \frac{1-F(M)}{F(M)} \mathbb{E}\#\{e \in \underline{\textsc{GEO}}(0,x) : D_e > M\}. \label{eq: D_e_inequality}
\end{align}

For the rest of the proof, we estimate the expected value in \eqref{eq: D_e_inequality}. Suppose first that 
\[
\text{the distribution of }t_e\text{ is bounded}.
\]
Then by \eqref{eq: precise_geodesic_moment_bound}, for some $\Cl[lgc]{c: android_constant}>0$, we have
\[
\mathbb{E}\#\{e \in \underline{\textsc{GEO}}(0,x) : D_e > M\} \leq \mathbb{E}\#\textsc{GEO}(0,x) \leq \Cr{c: android_constant}|x|,
\]
and so
\[
\mathbb{E}\overline{N}(0,x,(M,\infty)) \leq \Cr{c: android_constant} \frac{1-F(M)}{F(M)} |x|.
\]
This implies the statement of Thm.~\ref{thm: upper_bound_upper_tail} for all $M$ such that $F(M) \in (0,1)$. (Here we use that for fixed $\Cr{c: L2},\Cr{c: L1}>0$, the quantity $q(M,x)/|x|$ is bounded away from zero when $F(M) \in (0,1)$ because the set of $M$ for which $F(M) \in (0,1)$ is bounded.)

We may therefore suppose that
\begin{equation}\label{eq: distribution_unbounded}
\text{the distribution of }t_e \text{ is unbounded}.
\end{equation}
In this case, we will show that there exist 
$\Cr{c: L2}, \Cr{c: L1}>0$
such that for all $x \in \mathbb{Z}^d$ and $M$ with $F(M) \in (0,1)$, we have
\begin{equation}\label{eq: overall_inequality}
\mathbb{E}\#\{e \in \underline{\textsc{GEO}}(0,x) : D_e > M\} \leq q(M,x).
\end{equation}
Once we prove this inequality, we can combine it with \eqref{eq: D_e_inequality} to complete the proof. 

To prove \eqref{eq: overall_inequality}, we focus on large $M$, since we will use \eqref{eq: precise_geodesic_moment_bound} for small $M$.
The main idea of the proof of \eqref{eq: overall_inequality} for large $M$ is that if an edge $e$ is in a geodesic $\gamma$ from $0$ to $x$ but $D_e > M$, then we cannot find a path $\pi$ that connects two points of $\gamma$ on either side of $e$ (that is, $\pi$ ``shortcuts'' $e$) but has $T(\pi) \leq M$. We will show that the probability of existence of such a detour is large when $M$ is large, and so it is unlikely to find any path with too many edges that do not have shortcuts.

To begin, we make precise what we mean by existence of shortcuts. For an edge $e = \{u,v\}$ and $k \geq 1$, write $B(e,k)$ for the set of vertices $y$ with $\|y - (1/2)(u+v)\|_\infty \leq k$. Define the external vertex boundary by 
\begin{equation}\label{eq: boundary_definition}
\partial B(e,k) = \{y \in \mathbb{Z}^d \setminus B(e,k): \exists~z \in B(e,k) \text{ with } |z-y|=1\}.
\end{equation} 

\begin{figure}[h]
  \centering
  \includegraphics[width=0.6\linewidth, trim={5cm 12cm 5cm 6cm}, clip]{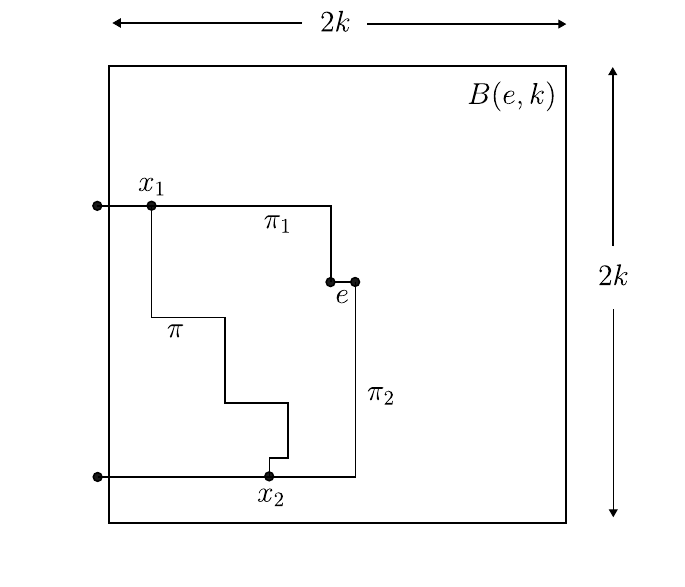}
  \caption{Illustration of the event that the edge $e$ is good. The paths $\pi_1$ and $\pi_2$ begin at $e$ and end at vertices in $\partial B(e,k)$. The path $\pi$ connects a vertex $x_1$ on $\pi_1$ to a vertex $x_2$ on $\pi_2$, and has at most $\Cr{c: good_edge} k$ many edges. All edges on $\pi$ have weight at most $M_0$ and $e$ is not an edge of $\pi$.}
  \label{fig: fig_1}
\end{figure}

\begin{df}\label{def: good}
Given numbers $\Cl[lgc]{c: good_edge},M_0 > 0$ and an integer $k \geq 1$, the edge $e$ is said to be \underline{good} if for each pair of paths $\pi_1,\pi_2$ in $B(e,k)$, starting at endpoints of $e$ and ending in the set $\partial B(e,k)$, there exist vertices $x_1 \in \pi_1$ and $x_2 \in \pi_2$ with the following property. There is a path $\pi$ from $x_1$ to $x_2$ such that
\begin{enumerate}
\item $e \notin \pi$,
\item $\pi$ has at most $\Cr{c: good_edge} k$ many edges, and
\item $t_f \leq M_0$ for all edges $f$ of $\pi$.
\end{enumerate}
\end{df}
See Figure~\ref{fig: fig_1}. We stress that the path $\pi$ from Def.~\ref{def: good} need not lie in $B(e,k)$. Because $e \notin \pi$, the event $\{e \text{ is good}\}$ is independent of $t_e$. We now estimate the probability that an edge is good.
\begin{lem}\label{lem: good_lemma}
The following hold for all edges $e$.
\begin{enumerate}
\item For $k=1, \Cr{c: good_edge}=3,$ and all $M_0>0$, we have
\[
\mathbb{P}(e \text{ is good}) \geq 1 - (3(1-F(M_0)))^{2d-2}.
\]
\item There exist $\Cr{c: good_edge} \geq 1$ and $\Cl[smc]{c: good_lemma},\Cl[lgc]{c: new_constant_2}, M_0>0$ such that for all $k \geq 1$, we have
\[
\mathbb{P}(e \text{ is good}) \geq 1 - \Cr{c: new_constant_2} e^{-\Cr{c: good_lemma}k}.
\]
\end{enumerate}
\end{lem}
\begin{proof}
We start with item 1, and we may take $e$ to have endpoints $0$ and $\mathbf{e}_1 = (1, 0, \dots, 0)$. Enumerate the neighbors of 0 (the vertices $\pm \mathbf{e}_i$) as $y_1, \dots, y_{2d}$ in any way that $y_1 = \mathbf{e}_1$ and $y_2 = -\mathbf{e}_1$. Define (edge-disjoint) paths $\gamma_3, \dots, \gamma_{2d}$ from 0 to $\mathbf{e}_1$ so that $\gamma_i$ starts at 0, moves to $y_i$, then to $y_i + \mathbf{e}_1$, and then to $\mathbf{e}_1$. Observe that if any path $\gamma_i$ has all edges $f$ with $t_f \leq M_0$, then $e$ is good. (Take $x_1=0,x_2=\mathbf{e}_1$ in Def.~\ref{def: good}.) Using independence, we obtain
\[
\mathbb{P}(e \text{ is not good}) \leq \prod_{i=3}^{2d} \mathbb{P}(\text{some }f \text{ in } \gamma_i \text{ has } t_f > M_0) \leq (3(1-F(M_0)))^{2d-2}.
\]

We now move to item 2. As we did in the proof of Prop.~\ref{prop: improved_tail_bound}, letting $M_0$ be such that $F(M_0) \geq p_0$ (from Lem.~\ref{lem: gamma_lemma}), we call an edge $e$ \underline{open} if $t_e \leq M_0$; otherwise, it is \underline{closed}. We again let $I$ be the unique infinite open cluster and define $A_{e,k}$ to be the event 
\[
A_{e,k} = \{\text{every path from }e\text{ to }\partial B(e,k)\text{ intersects }I\}.
\]
By Lem.~\ref{lem: gamma_lemma}, we have 
\begin{equation}\label{eq: normal_bound}
\mathbb{P}(A_{e,k}) \geq 1 - e^{-\Cr{c: book_chemical_constant}k} \text{ for all } k \geq 1.
\end{equation}
Using \eqref{eq: distribution_unbounded}, we henceforth fix $M_0$ such that
\[
F(M_0) \in [p_0,1).
\]

In addition to the event $A_{e,k}$, we will use $A_{e,k}'$, defined as
\[
A_{e,k}' = \{\forall x,y \in I \cap B(e,2k), \text{ we have } d_I(x,y) \leq \Cl[lgc]{c: chemical}k\}.
\]
Here, $\Cr{c: chemical}>0$ is a constant and, as before, $d_I(x,y)$ is the minimal number of edges in any open path from $x$ to $y$, with $d_I(x,y)= \infty$ if there is no such path. To estimate the probability of $A'_{e,k}$, we directly apply Lem.~\ref{lem: AP_lem}. Choose $\Cr{c: chemical}$ from the definition of $A_{e,k}'$ to be such that $\Cr{c: chemical} k \geq \Cr{c: AP_large_constant}|x-y|$ for all $x,y \in B(e,2k)$, where $\Cr{c: AP_large_constant}$ is from Lem.~\ref{lem: AP_lem}. Then for $x,y \in B(e,2k)$, we have $\mathbb{P}( \Cr{c: chemical}k \leq d_I(x,y) < \infty)\leq e^{-\Cr{c: AP_small_constant} \Cr{c: chemical} k}$, and so 
we find for constants $\Cl[smc]{c: L22},\Cl[lgc]{c: L20},\Cl[lgc]{c: L21}>0$,
\begin{equation}\label{eq: prime_bound}
\mathbb{P}((A_{e,k}')^c) \leq \sum_{x,y \in B(e,2k)} \mathbb{P}(\Cr{c: chemical} k \leq d_I(x,y) < \infty) \leq \Cr{c: L20} k^{2d} e^{-\Cr{c: AP_small_constant} \Cr{c: chemical} k} \leq \Cr{c: L21}e^{- \Cr{c: L22} k}.
\end{equation}

Let us say that $e$ is \underline{nearly good} if all of the conditions of Def.~\ref{def: good} hold, except possibly item 1. That is, the path $\pi$ is permitted to contain $e$. We claim that
\begin{equation}\label{eq: detour_claim}
A_{e,k} \cap A_{e,k}' \subset \{e \text{ is nearly good}\}.
\end{equation}
To prove this, consider an outcome in $A_{e,k} \cap A_{e,k}'$. Let $\pi_1,\pi_2$ be paths in $B(e,k)$ starting at endpoints of $e$ and ending at $\partial B(e,k)$. Since $A_{e,k}$ occurs, there exist vertices $x_1 \in \pi_1$ and $x_2 \in \pi_2$ such that $x_1,x_2 \in I$. Both of $x_1$ and $x_2$ are in $B(e,2k)$, so since $A_{e,k}'$ occurs, we have $d_I(x_1,x_2) \leq \Cr{c: chemical}k$. This shows that $e$ is nearly good, so long as we choose $\Cr{c: good_edge}$ from Def.~\ref{def: good} to equal $\Cr{c: chemical}$.

Last, we observe that
\begin{equation}\label{eq: labelhead}
\{e \text{ is not good}\}\cap \{t_e > M_0\} \subset \{e \text{ is not nearly good}\}.
\end{equation}
Indeed, if $t_e > M_0$ but $e$ is nearly good, then any path $\pi$ in the definition of nearly good cannot contain $e$, and therefore $e$ is good. By independence, \eqref{eq: normal_bound}, \eqref{eq: prime_bound}, and \eqref{eq: detour_claim},
\begin{align*}
\mathbb{P}(e \text{ is not good}) &= \frac{1}{1-F(M_0)} \mathbb{P}(e \text{ is not good}, t_e > M_0) \\
&\leq \frac{1}{1-F(M_0)} \mathbb{P}(e \text{ is not nearly good}) \\
&\leq \frac{1}{1-F(M_0)}\left( \mathbb{P}(A_{e,k}^c) + \mathbb{P}((A_{e,k}')^c) \right) \\
&\leq \frac{1}{1-F(M_0)}\left( e^{-\Cr{c: book_chemical_constant} k} + \Cr{c: L21} e^{-\Cr{c: L22} k} \right).
\end{align*}
This completes the proof of the lemma.
\end{proof}

To relate the definition of good to the condition $D_e > M$ in \eqref{eq: D_e_inequality}, we give the following claim. Let $\Cr{c: good_edge}, M_0$ be as in Def.~\ref{def: good}.
\begin{claim}\label{claim: good_claim}
Let $e$ be an edge and $x \in \mathbb{Z}^d$. If $k \geq 1$ and $M$ satisfy
\begin{enumerate}
\item $\Cr{c: good_edge}kM_0 < M \leq D_e,$ and
\item $B(e,k)$ does not contain $0$ or $x$,
\end{enumerate}
then $e$ is not good.
\end{claim}
\begin{proof}
Consider a configuration $(t_f)$ for which $D_e \geq M$, and suppose that items 1 and 2 hold. Define the configuration $(\bar{t}_f)$ by $\bar{t}_f = M$ if $f=e$ and $\bar{t}_f = t_f$ if $f \neq e$. Because $D_e$ depends only on $t_e^c$, we have $D_e \geq M$ in $(\bar{t}_f)$. From \eqref{eq: overline_equivalence}, $e \in \overline{\textsc{GEO}}(0,x)$ in $(\bar{t}_f)$. Thus we may find a path $\gamma$ from $0$ to $x$ that contains $e$ and is a geodesic in $(\bar{t}_f)$. Since $0$ and $x$ are not in $B(e,k)$, $\gamma$ must intersect $\partial B(e,k)$ before and after touching $e$. Let $\pi_1$ be the segment of $\gamma$ from its last intersection with $\partial B(e,k)$ before it touches $e$ until it meets $e$ at an endpoint, and let $\pi_2$ be the segment of $\gamma$ from the other endpoint of $e$ until its first intersection with $\partial B(e,k)$. 

Assume for a contradiction that $e$ is good (in the original configuration $(t_f)$). Then we can pick $x_1 \in \pi_1$, $x_2 \in \pi_2$, and a path $\pi$ from $x_1$ to $x_2$ such that $e \notin \pi$, $\pi$ has at most $\Cr{c: good_edge} k$ many edges, and $t_f \leq M_0$ for all edges $f$ of $\pi$. Because $t_f = \bar{t}_f$ for all $f \neq e$, we have $\bar{t}_f \leq M_0$ for all $f$ in $\pi$. Writing $\bar{T}(x_1,x_2)$ for the passage time in the configuration $(\bar{t}_f)$, we have 
\[
\bar{T}(x_1,x_2) \geq \bar{t}_e = M,
\]
but $\bar{T}(x_1,x_2) \leq \bar{T}(\pi) \leq M_0\Cr{c: good_edge} k < M$. This contradiction proves that $e$ is not good.
\end{proof}

We now return to proving \eqref{eq: overall_inequality}, and for this, we suppose for now that $k \geq 1, M>0,$ and $M_0>0$ are such that $M > \Cr{c: good_edge} k M_0$. 
Due to Claim~\ref{claim: good_claim}, we get for some $\Cl[lgc]{c: another_constant}, \Cl[lgc]{c: another_constant_2}>0$,
\begin{align}
&\mathbb{E}\#\{e \in \underline{\textsc{GEO}}(0,x) : D_e>M\} \nonumber \\
\leq~& \mathbb{E}\#\{e \in \textsc{GEO}(0,x) : D_e > M\} \nonumber \\
\leq~& \#\{e : 0 \in B(e,k) \text{ or } x \in B(e,k)\} \nonumber \\
+~& \mathbb{E}\#\{e \in \textsc{GEO}(0,x) : e \text{ is not good}\} \nonumber \\
\leq~& \Cr{c: another_constant} k^d +  \mathbb{E}\#\{e \in \textsc{GEO}(0,x) : e \text{ is not good}\} \nonumber \\
\leq~& \Cr{c: another_constant} k^d + \Cr{c: another_constant_2} +  \mathbb{E}\#\{e \in \textsc{GEO}(0,x) : e \text{ is not good}\}\mathbf{1}_{\{\#\textsc{GEO}(0,x) \leq \Cr{c: improved_big_constant} |x|\}} \label{eq: pasta_macaroni}.
\end{align}
In the final inequality we have used the following consequence of Prop.~\ref{prop: improved_tail_bound}: for some $\Cr{c: another_constant_2}>0$ and all nonzero $x$, we have
\[
\mathbb{E}\#\textsc{GEO}(0,x) \mathbf{1}_{\{\#\textsc{GEO}(0,x) \geq \Cr{c: improved_big_constant}|x|\}} \leq \Cr{c: another_constant_2}.
\]

To estimate \eqref{eq: pasta_macaroni}, we use the lattice animal bound from Appendix~\ref{sec: appendix}. For $n \geq 1$, let $\Xi_n$ be the collection of subsets of $\mathcal{E}^d$ of cardinality $n$ which contain an edge incident to the origin and which are connected (a set $B$ of edges is connected if for each $e,f \in B$, there is a path with initial edge $e$ and final edge $f$, all whose edges lie in $B$). Define
\[
N_n = \max_{S \in \Xi_n} \left( \sum_{e \in S} \mathbf{1}_{\{e \text{ is not good}\}}\right).
\]
Continuing from \eqref{eq: pasta_macaroni}, we find for $M > \Cr{c: good_edge}k M_0$ and some $\Cl[lgc]{c: more_constants}>0$,
\begin{equation}\label{eq: new_taco_head}
\mathbb{E}\#\{e \in \underline{\textsc{GEO}}(0,x) : D_e>M\} \leq \Cr{c: more_constants} k^d + \mathbb{E}N_{\lfloor \Cr{c: improved_big_constant} |x|\rfloor}.
\end{equation}
The variables $(\mathbf{1}_{\{e \text{ is not good}\}})$ are $2(1+\Cr{c: good_edge})k$-dependent in the sense of Definition~\ref{def: k_dependent}. Taking $p$ as the probability that an edge is not good, we apply Prop.~\ref{prop: lattice_animals} to obtain
\[
 \mathbb{E}N_{\lfloor \Cr{c: improved_big_constant} \|x\|_1\rfloor} \leq \Cr{c: lattice_animal_constant} (2(1+\Cr{c: good_edge})k)^d p^{\frac{1}{d}} \Cr{c: improved_big_constant} |x|. 
 \]
We put this in \eqref{eq: new_taco_head} to get $\Cl[lgc]{c: near_the_end}>0$ such that for $M > \Cr{c: good_edge}k M_0$,
\begin{equation}\label{eq: to_apply_in_the_end}
\mathbb{E}\#\{e \in \underline{\textsc{GEO}}(0,x) : D_e>M\} \leq \Cr{c: near_the_end}k^d \left( 1+ p^{\frac{1}{d}} |x|\right).
\end{equation}

If we apply \eqref{eq: to_apply_in_the_end} with $k=1, \Cr{c: good_edge}=3$, $M>0$, and $M_0 = M/4$, then item 1 of Lem.~\ref{lem: good_lemma} gives $p \leq (3(1-F(M/4)))^{2d-2}$, and so
\begin{equation}\label{eq: first_end_bound}
\mathbb{E}\#\{e \in \underline{\textsc{GEO}}(0,x) : D_e>M\} \leq \Cr{c: near_the_end} \left( 1+ (3(1-F(M/4)))^{2- \frac{2}{d}} |x|\right) \text{ for }M>0.
\end{equation}
Next, taking $\Cr{c: good_edge}$ and $M_0$ from item 2 of Lem.~\ref{lem: good_lemma}, we now suppose that
\[
M \geq 2M_0 \Cr{c: good_edge}
\]
and pick 
\[
k = \lfloor M/(2M_0\Cr{c: good_edge}) \rfloor \geq 1,
\]
so that $M_0\Cr{c: good_edge} k < M$. Then from that item, we have $p \leq \Cr{c: new_constant_2} e^{-\Cr{c: good_lemma}k}$, and so from \eqref{eq: to_apply_in_the_end}, we obtain for some $\Cl[smc]{c: david_delight_2},\Cl[lgc]{c: david_delight_1}>0$,
\begin{equation}\label{eq: second_end_bound}
\mathbb{E}\#\{e \in \underline{\textsc{GEO}}(0,x) : D_e>M\} \leq \Cr{c: david_delight_1} M^d \left( 1+ e^{-\Cr{c: david_delight_2}M} |x|\right) \text{ for } M \geq 2M_0 \Cr{c: good_edge}.
\end{equation}
Last, for $M \leq 2M_0 \Cr{c: good_edge}$, we instead use \eqref{eq: precise_geodesic_moment_bound} to obtain $\Cl[lgc]{c: constant_tempura}>0$ such that for all nonzero $x$, we have 
\[
\mathbb{E}\#\{e \in \underline{\textsc{GEO}}(0,x) : D_e>M\} \leq \mathbb{E}\#\textsc{GEO}(0,x) \leq \Cr{c: constant_tempura}|x|.
\]
Combining this with \eqref{eq: first_end_bound} and \eqref{eq: second_end_bound} gives \eqref{eq: overall_inequality} and completes the proof.

\section{Proof of Thm.~\ref{thm: BK_argument}}\label{sec: BK_argument}

The proof will be split into steps.

\medskip
\noindent
{\bf Step 1: Basic definitions and black rectangles.} To show that a geodesic has many edges with large weight $\geq M$, we will consider crossings of rectangles made by a geodesic, and modify the weights in these rectangles to force the geodesic to take high-weight edges. We start with some definitions.

Let $n \geq 1$ and define the rectangular box
\[
R(n) = [0,n] \times [0,3n]^{d-1}.
\]
We say that a vertex self-avoiding path $\pi$ \underline{is a crossing of $R(n)$} if, when we write the vertices of $\pi$ in order as $x_0, \dots, x_k$, we have $x_\ell \in R(n)$ for all $\ell$, $x_0 \cdot \mathbf{e}_1 = 0$, and $x_k \cdot \mathbf{e}_1 = n$. We say that $\pi$ \underline{crosses} $R(n)$ if it contains a subpath which is a crossing of $R(n)$.

Our weight modification will require the rectangle $R(n)$ to satisfy two properties, both of which hold with high probability when $n$ is large. To state them, we need a few more definitions. We say that a vertex self-avoiding path $\pi$, with vertices $x_0, \dots, x_k$ in order, is \underline{oriented} if $k = \|x_0-x_k\|_1$, and we say that $\pi$ is \underline{mostly oriented} if $\pi$ has an oriented subpath (that is, an oriented path that begins at some $x_j$ and follows $\pi$ until it stops at some $x_{j'}$ with $j' \geq j$) with at least $k/3$ many edges. We say that $\pi$ is a \underline{rung} for $R(n)$ if
\begin{enumerate}
\item $\pi$ is mostly oriented,
\item $x_\ell \in R(n)$ for all $\ell$, and 
\item $x_\ell \in \partial_i R(n)$ if and only if $\ell = 0$ or $k$, where $\partial_iR(n)$ is the internal vertex boundary $\{w \in R(n) : \exists v \in R(n)^c \text{ such that } |w-v|=1\}$.
\end{enumerate}

We now list the properties $R(n)$ must satisfy. Recall that $r$ is the infimum of the support of the distribution of $t_e$.
\begin{df}\label{def: black}
Let $\delta, \Cl[lgc]{c: black_constant}>0$. We say that the rectangle $R(n)$ is \underline{black} if 
\begin{enumerate}
\item there is no rung for $R(n)$, with at least $\Cr{c: black_constant} \log n$ many edges, all of whose edges $e$ satisfy $t_e \leq r+\delta/2$, and
\item for all vertices $a,b \in R(n)$ with $\|a-b\|_1 \geq \Cr{c: black_constant} \log n$, we have $T(a,b) \geq (r+\delta)\|a-b\|_1$.
\end{enumerate}
\end{df}

Choose any $r'$ such that
\begin{equation}\label{eq: r_prime}
r' > r \text{ and } \mathbb{P}(t_e \in (r,r']) > 0.
\end{equation}
Regarding the probability of a rectangle being black, we have the following lemma.

\begin{lem}\label{lem: black}
There exist $\delta,\Cr{c: black_constant}>0$  
such that 
\[
\mathbb{P}(t_e \in [r+\delta,r'])>0
\]
and
\[
\mathbb{P}(R(n) \text{ is black}) \to 1 \text{ as } n \to \infty.
\]
\end{lem}
\begin{proof}
Because $\mathbb{P}(t_e \in (r,r'])>0$, we can choose $\delta_1>0$ such that $\mathbb{P}(t_e \in [r+\delta_1,r'])>0$. To show that $\mathbb{P}(R(n) \text{ is black}) \to 1$, we first use \cite[Lem.~5.3]{BK93}, which states that under assumption \eqref{eq: useful}, there exist $\delta_2>0$ and $\Cl[smc]{c: BK_constant_1}>0$ such that for all vertices $v,w$, we have
\[
\mathbb{P}(T(v,w) < (r+\delta_2) \|v-w\|_1) \leq e^{-\Cr{c: BK_constant_1} \|v-w\|_1}.
\]
This with a union bound implies
\begin{align*}
\mathbb{P}(\text{item }2 \text{ fails for }\delta = \delta_2) &\leq \sum_{a,b \in R(n) : \|a-b\|_1 \geq \Cr{c: black_constant} \log n} \mathbb{P}(T(a,b) < (r+\delta_2) \|a-b\|_1) \\
&\leq \sum_{a,b \in R(n) : \|a-b\|_1 \geq \Cr{c: black_constant} \log n} e^{-\Cr{c: BK_constant_1} \|a-b\|_1} \\
&\leq (3n)^{2d} \exp\left( - \Cr{c: black_constant} \Cr{c: BK_constant_1} \log n\right) \to 0 \text{ as } n \to \infty,
\end{align*}
so long as we choose $\Cr{c: black_constant} > 2d/\Cr{c: BK_constant_1}$.

Turning to item 1, if there is a rung for $R(n)$ with at least $\Cr{c: black_constant} \log n$ many edges, all whose edges $e$ satisfy $t_e \leq r+\delta_2/2$, then because rungs are mostly oriented, there must exist a vertex $v \in R(n)$ from which there emanates an oriented path $\pi$ in $R(n)$ with $\#\pi \geq (\Cr{c: black_constant} \log n) / 3$, and all edges $e$ of $\pi$ satisfy $t_e \leq r + \delta_2/2$. If $w \in R(n)$ is the other endpoint of $\pi$, then $\|v-w\|_1 \geq (\Cr{c: black_constant} \log n)/3$ and $T(v,w) \leq T(\pi) \leq (r+\delta_2/2)\|v-w\|_1$. So by the same argument as above, we have
\begin{align*}
\mathbb{P}(\text{item }1 \text{ fails for } \delta = \delta_2) &\leq \sum_{v,w \in R(n), \|v-w\|_1 \geq (\Cr{c: black_constant} \log n)/3} \mathbb{P}\left( T(v,w) \leq \left( r+ \frac{\delta_2}{2}\right) \|v-w\|_1\right) \\
&\leq (3n)^{2d} \exp\left( - \left( \Cr{c: black_constant} \Cr{c: BK_constant_1}/3\right) \log n\right) \to 0 \text{ as } n \to \infty,
\end{align*}
so long as we choose $\Cr{c: black_constant} > 6d/\Cr{c: BK_constant_1}$. Taking $\delta = \min \{ \delta_1, \delta_2\}$ completes the proof.
\end{proof}

\medskip
\noindent
{\bf Step 2: Edge-modification argument.}
We henceforth fix the values of $\delta$ and $\Cr{c: black_constant}$ from Lem.~\ref{lem: black}. Given vertices $x,y \in \mathbb{Z}^d$, let $\gamma_{x,y}$ be the vertex self-avoiding geodesic between $x$ and $y$. If there are multiple such geodesics, choose the first in some deterministic ordering of all finite paths in $\mathbb{Z}^d$. This $\gamma_{x,y}$ is defined almost surely under assumption~\eqref{eq: percolation_condition}, which is implied by \eqref{eq: useful}. For $n \geq 1$, let $\overline{R}(n)$ be the set of edges which have at least one endpoint in $R(n)$. We will now argue that if $R(n)$ is black and $x,y$ are vertices not in $R(n)$ such that $\gamma_{x,y}$ crosses $R(n)$, then we can modify the weights in $\overline{R}(n)$ so as to force $\gamma_{x,y}$ to take an edge in $\overline{R}(n)$ with weight $\geq M$.

\begin{prop}\label{prop: modification}
There exists $\Cl[lgc]{c: modification_constant}>0$ such that the following holds. For any $x,y \in \mathbb{Z}^d \setminus R(n)$ and $n,M$ satisfying $n > 32M/\delta$ with $M$ sufficiently large, we have
\begin{align*}
&\mathbb{P}(\text{each geodesic from } x \text{ to }y \text{ contains an edge } e \in \overline{R}(n) \text{ with } t_e \geq M) \\
\geq~&\mathbb{P}(t_e \in [M,2M])^{2d} \mathbb{P}(t_e \geq M)^{\Cr{c: modification_constant} n} \mathbb{P}(\gamma_{x,y} \text{ crosses }R(n) \text{ and } R(n) \text{ is black}).
\end{align*}
\end{prop}
\begin{proof}
On the event $\{\gamma_{x,y} \text{ crosses } R(n) \text{ and } R(n) \text{ is black}\}$, define $u_0$ as the first vertex that $\gamma_{x,y}$ takes in $R(n)$ (as we proceed from $x$ to $y$ along $\gamma_{x,y}$), and define $v_0$ as the last vertex that $\gamma_{x,y}$ takes in $R(n)$. We decompose the probability as
\begin{align}
&\mathbb{P}(\gamma_{x,y} \text{ crosses } R(n),~R(n) \text{ is black}) \nonumber  \\
=~& \sum_{u,v \in \partial_iR(n)} \mathbb{P}\left(u_0 = u, v_0 = v,~ \gamma_{x,y}\text{ crosses } R(n),~R(n) \text{ is black}\right). \label{eq: my_decomposition}
\end{align}
We aim to replace the portion of the path $\gamma_{x,y}$ between $u$ and $v$ with a deterministic rung $\pi_{u,v}$. However, as a rung can only touch $\partial_i R(n)$ twice, this will be impossible if $u$ or $v$ lie on certain parts of the boundary $\partial_i R(n)$ (for example, if $d=2$ and $u=(0,3n)$). Therefore we must choose new endpoints $u^\ast$ and $v^\ast$ of the rung. For this purpose, choose deterministic (that is, dependent only on $u$ and $v$) vertices $u^\ast,v^\ast \in \partial_iR(n)$ satisfying: 
\begin{itemize}
\item there are vertex disjoint paths $\pi_{u,v}^{(1)}$ connecting $u$ to $u^\ast$ and $\pi_{u,v}^{(2)}$ connecting $v^\ast$ to $v$, remaining in $\partial_iR(n)$ and each having at most $d-1$ many edges, and
\item $u^\ast$ is adjacent to a vertex $\overline{u} \in R(n) \setminus \partial_iR(n)$ and $v^\ast$ is adjacent to a different vertex $\overline{v} \neq \overline{u}$ in $R(n) \setminus \partial_iR(n)$.
\end{itemize}
For example, if $u = (n,0)$ and $v=(n,1)$, we can choose $u^\ast = (n-1,0), \overline{u} = (n-1,1)$ and $v^\ast = (n,2), \overline{v} = (n-1,2)$, and $\pi^{(i)}_{u,v}$ for $i=1,2$ to be paths consisting of one edge each. See the top two subfigures of Fig.~\ref{fig: fig_2} for an illustration of these and the following definitions. There, we have chosen $u = (0,3n), u^\ast = (1,3n), \overline{u} = (1,3n-1)$, and $v = v^\ast$. 

\begin{figure}[p]
  \hspace{-.7in}
  \begin{subfigure}[b]{0.4\linewidth}
  \includegraphics[width=2\linewidth,trim={2cm 18cm 0cm .75cm}, clip]{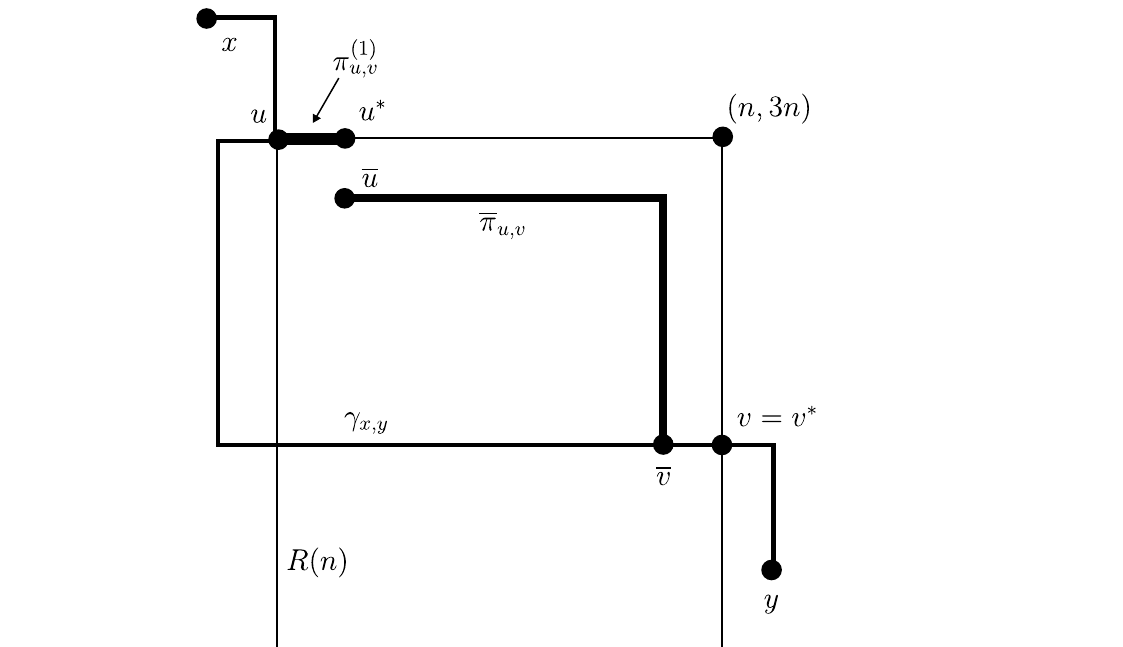}
  \end{subfigure}
  \begin{subfigure}[b]{0.4\linewidth}
  \includegraphics[width=2\linewidth,trim={0cm 18cm 2cm .75cm}, clip]{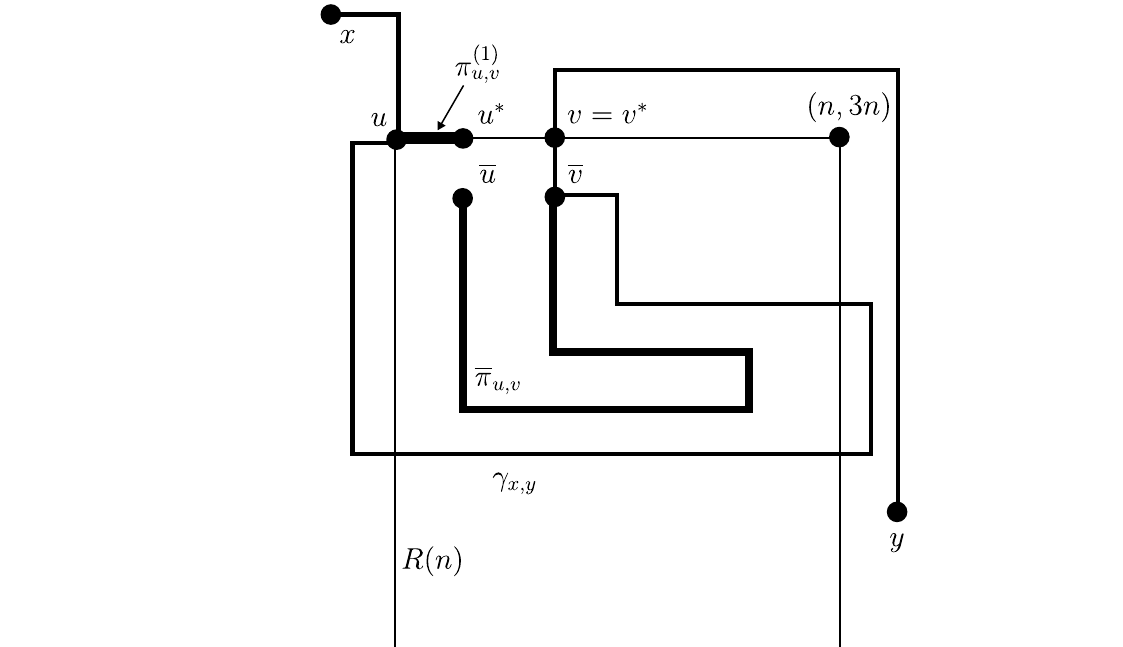}
  \end{subfigure}
  
  \medskip
  \hspace{-.7in}
  \begin{subfigure}[b]{0.4\linewidth}
  \includegraphics[width=2\linewidth,trim={2cm 15.5cm 0cm 1cm}, clip]{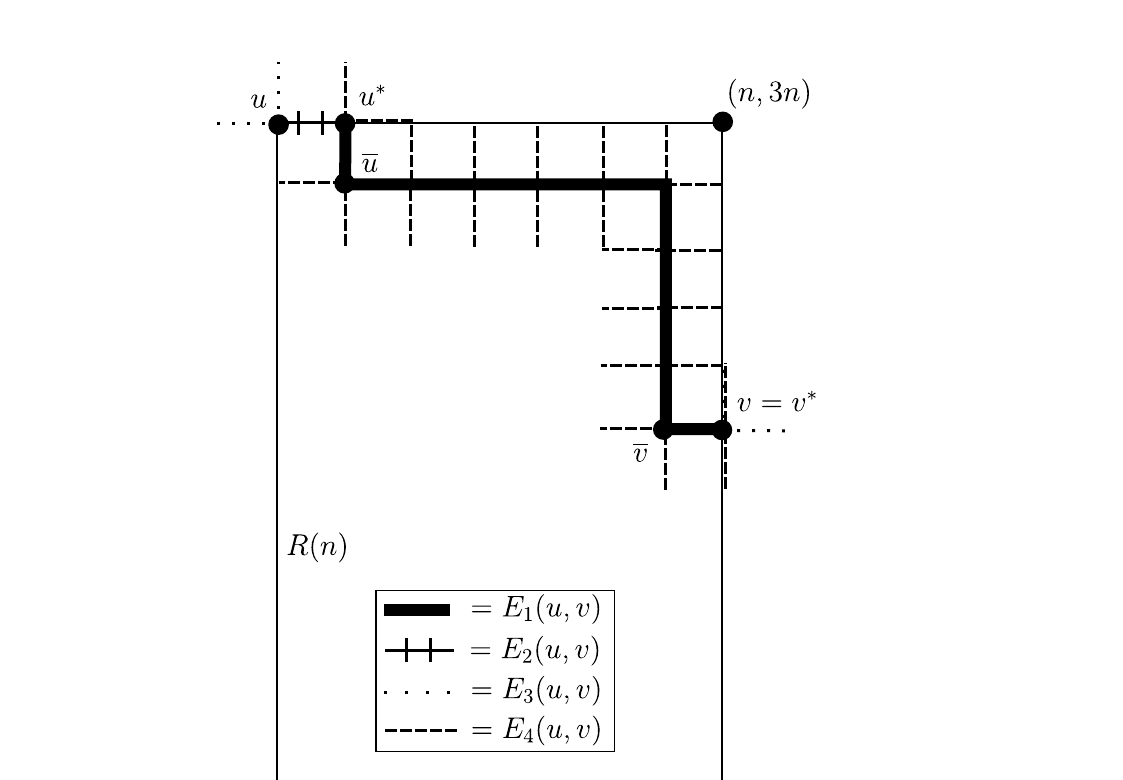}
  \end{subfigure}
  \begin{subfigure}[b]{0.4\linewidth}
  \includegraphics[width=2\linewidth,trim={0cm 15.5cm 2cm 1cm}, clip]{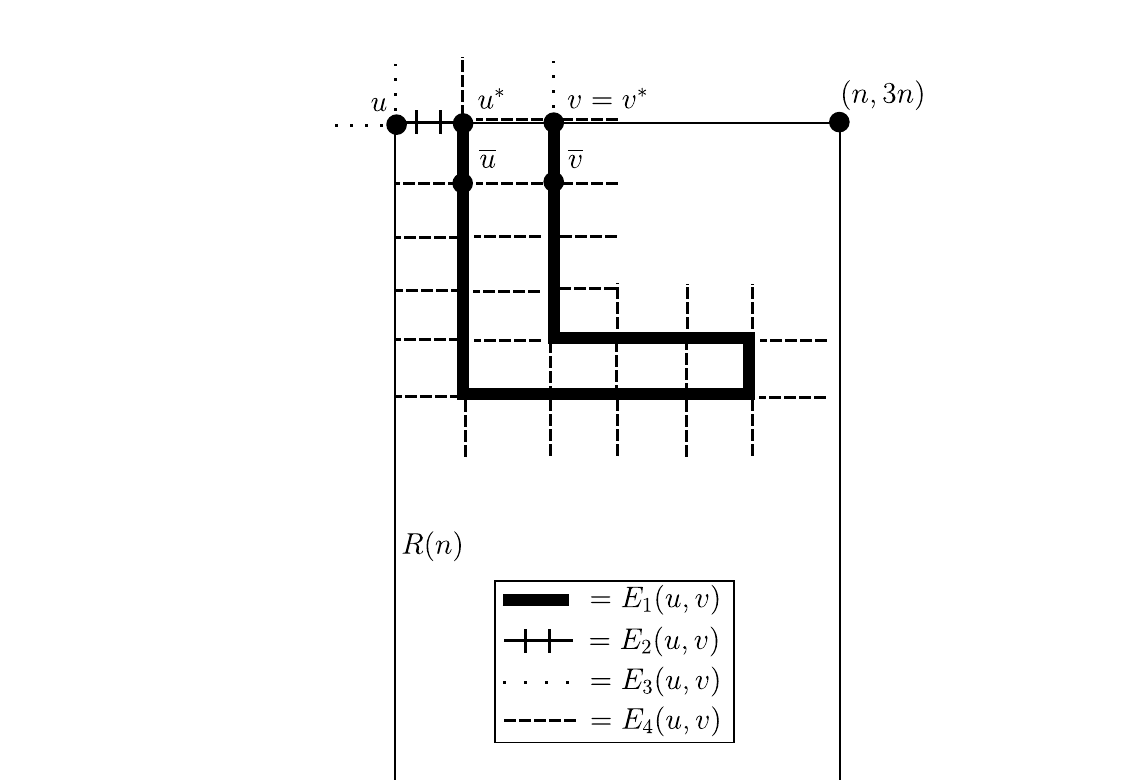}
  \end{subfigure}
  \caption{Illustration of path and edge set definitions. In top-left, $u$ and $v$, the first and last intersections of $\gamma_{x,y}$ with $R(n)$, are far apart ($\|u-v\|_1 \geq n/2$) and so the path $\overline{\pi}_{u,v}$ from $\overline{u}$ to $\overline{v}$ is chosen to be oriented. In top-right, $u$ and $v$ are close ($\|u-v\|_1 < n/2$) and so $\overline{\pi}_{u,v}$ is a concatenation of two oriented paths. In both cases, $u$ is the corner vertex $(0,3n)$ and is not adjacent to $R(n) \setminus \partial_iR(n)$, so it cannot be the endpoint of a rung, and we choose $u^\ast = (1,3n)$. In both cases, we choose $v = v^\ast$ and so $\pi_{u,v}^{(2)}$ is empty. The rung $\pi_{u,v}$ (not labeled) moves from $u^\ast$ to $\overline{u}$, follows $\overline{\pi}_{u,v}$, then moves from $\overline{v}$ to $v^\ast$. The path $\sigma_{u,v}$ (not labeled) follows $\pi_{u,v}^{(1)}$, then $\pi_{u,v}$, and then $\pi_{u,v}^{(2)}$. The bottom two subfigures illustrate the sets $E_i(u,v)$ for $i=1, \dots, 4$ and correspond to the cases shown in the subfigures directly above them. Note that $\overline{\pi}_{u,v}$ has length $<n$ in the two right figures, although due to space constraints, the path appears longer than the width $n$ of the box.}
  \label{fig: fig_2}
\end{figure}

We now define the rung $\pi_{u,v}$; choose it to be any deterministic vertex self-avoiding path satisfying the following conditions:
\begin{itemize}
\item the first edge of $\pi_{u,v}$ connects $u^\ast$ to $\overline{u}$ and the last edge connects $\overline{v}$ to $v^\ast$,
\item the subpath $\overline{\pi}_{u,v}$ of $\pi_{u,v}$ between $\overline{u}$ and $\overline{v}$ is contained in $R(n) \setminus \partial_i R(n)$,
\item if $\|u-v\|_1 \geq n/2$, then $\overline{\pi}_{u,v}$ is oriented, and
\item if $\|u-v\|_1 < n/2$, then $\overline{\pi}_{u,v}$ consists of an oriented path followed by another oriented path, and has total length in the interval $[n/2, n)$.
\end{itemize}

We would like to choose the path $\overline{\pi}_{u,v}$ to be oriented to minimize its length and therefore minimize its passage time in our modified weights. If we were to do this in the case that $\|u-v\|_1$ is small, the length of $\pi_{u,v}$ would also be small, say less than $\Cr{c: black_constant} \log n$, where $\Cr{c: black_constant}$ is from Def.~\ref{def: black}. But in our modified weights, we will need $\pi_{u,v}$ to be the only rung that has all edges with weight at most $r + \delta/2$ (see \eqref{eq: what_we_need_for_chris}, where we need the displayed events to be disjoint as $\pi$ varies). This will be guaranteed from Def.~\ref{def: black} only when $\pi_{u,v}$ has length $\geq \Cr{c: black_constant} \log n$. So we must split off the case that $\|u-v\|_1$ is small, in the fourth bullet above, and only require $\overline{\pi}_{u,v}$ to be mostly oriented.

Observe that $\pi_{u,v}$ contains only two vertices (its first and last) in $\partial_iR(n)$. Furthermore, it is mostly oriented. Indeed, if $\|u-v\|_1 \geq n/2$, then $\pi_{u,v}$ is oriented except for at most its initial and final edges. If $\|u-v\|_1 < n/2$, then at least one of the two oriented subpaths in the fourth bullet point has at least $(k-2)/2$ many edges, where $k$ is the number of edges of $\pi_{u,v}$. If $n$ is large, then this is at least $k/3$, giving that $\pi_{u,v}$ is mostly oriented. We conclude that $\pi_{u,v}$ is a rung.

Finally, define $\sigma_{u,v}$ to be the path from $u$ to $v$ obtained by starting at $u$, following $\pi_{u,v}^{(1)}$ to $u^\ast$, following $\pi_{u,v}$ from $u^\ast$ to $v^\ast$, and then following $\pi_{u,v}^{(2)}$ from $v^\ast$ to $v$. 

To perform the edge modification, we need to define some edge sets. Let
\begin{align*}
E_1(u,v) &= \{e \in \mathcal{E}^d : e \text{ is an edge of } \pi_{u,v}\}, \\
E_2(u,v) &= \{e \in \mathcal{E}^d : e \text{ is an edge of }\pi_{u,v}^{(1)} \text{ or } \pi_{u,v}^{(2)}\}, \\
E_3(u,v) &= \{e \in \overline{R}(n) \setminus R(n) : e \text{ connects }u \text{ or }v \text{ to }  R(n)^c\}, \\
E_4(u,v) &= \{e \in \mathcal{E}^d \setminus \cup_{i=1}^3 E_i(u,v) : e \text{ is incident to } E_1(u,v) \cup E_2(u,v)\}.
\end{align*}

The edges of $E_1\cup E_2$ comprise the path $\sigma_{u,v}$, and those of $E_3\cup E_4$ comprise the edge boundary of $\sigma_{u,v}$, with those of $E_3$ being specifically adjacent to $u$ or $v$ from ``outside'' $R(n)$. All edges of $\cup_{i=1}^4 E_i(u,v)$ are in $\overline{R}(n)$. (See the bottom two subfigures of Fig.~\ref{fig: fig_2}.) The edges in $E_1$ will have low weight, those in $E_2$ will have typical weight, and those in $E_3 \cup E_4$ will have high weight, with those in $E_3$ having weight that is not too high. We make these specifications by sampling independent weights. Let $(t_e^\ast)$ be a family of i.i.d.~random variables, distributed identically to $t_e$, but independent of the original weights $(t_e)$, and define weights $(t_e^{u,v})$ by
\[
t_e^{u,v} = \begin{cases}
t_e &\quad\text{if } e \notin \cup_{i=1}^4 E_i(u,v) \\
t_e^\ast &\quad\text{if } e \in \cup_{i=1}^4 E_i(u,v).
\end{cases}
\]
Also, let $A_{u,v}^\ast$ be the event that the following conditions hold:
\begin{enumerate}
\item $t_e^\ast \leq r+\delta/2$ for $e \in E_1(u,v)$,
\item $t_e^\ast \in [r+\delta, r']$ for $e \in E_2(u,v)$,
\item $t_e^\ast \in [M,2M]$ for $e \in E_3(u,v)$, and
\item $t_e^\ast \geq M$ for $e \in E_4(u,v)$.
\end{enumerate}
To estimate the probability of $A_{u,v}^\ast$, we observe that $\#E_1(u,v) = 2+ \#\overline{\pi}_{u,v} \leq 2 + 3dn$. Also $\#E_2(u,v) \leq 2(d-1)$ and $\#E_3(u,v) \leq 2d$. Last, $\#E_4(u,v) \leq (4d-2) (\#E_1(u,v) + \#E_2(u,v)) \leq (4d-2)(2+3dn+2(d-1))$. In total, we can find $\Cl[lgc]{c: A_bound_constant}$ such that for all $n$ and all $u,v \in \partial_iR(n)$,
\begin{equation}\label{eq: A_bound}
\mathbb{P}(A_{u,v}^\ast) \geq e^{-\Cr{c: A_bound_constant} n} \mathbb{P}(t_e \in [M,2M])^{2d} \mathbb{P}(t_e \geq M)^{\Cr{c: A_bound_constant}n}.
\end{equation}

We claim that
\begin{align}
&\{\text{in }(t_e), u_0=u,v_0=v, \text{ and } R(n) \text{ is black}\} \cap A_{u,v}^\ast \nonumber \\
\subset~& \left\{ \text{in }(t_e^{u,v}), \begin{array}{c}
\text{each geodesic from }x \text{ to } y \text{ contains an edge} \\
e \in \overline{R}(n) \text{ with }t_e^{u,v} \geq M \text{ and } \pi_{u,v} \text{ is the only} \\
\text{rung with at least } \frac{n}{2} \text{ many edges,} \\
\text{all whose edges }e \text{ have } t_e^{u,v} \leq r+\frac{\delta}{2}
\end{array}
\right\}. \label{eq: main_containment}
\end{align}
(When we write $\pi_{u,v}$ is the only rung, we mean both $\pi_{u,v}$ and its reversal---the path obtained by following $\pi_{u,v}$ in reverse order.) To prove this, consider a configuration $(t_e)$ for which $u_0=u,v_0=v$, and $R(n)$ is black, and a configuration $(t_e^\ast)$ for which $A_{u,v}^\ast$ occurs.  Let $\tau$ be any geodesic from $x$ to $y$ in the weights $(t_e^{u,v})$. We will first prove that 
\begin{equation}\label{eq: take_large_edge}
\tau \text{ contains an edge of }E_3(u,v) \cup E_4(u,v).
\end{equation}
Note that, by definition of $A_{u,v}^\ast$, this will imply that $\tau$ contains an edge $e$ with $t_e^{u,v} \geq M$.

Let $\Gamma$ be the path $\gamma_{x,y}$ in the weights $(t_e)$ (which by assumption crosses $R(n)$), and define the path $\tau'$ from $x$ to $y$ that follows $\Gamma$ from $x$ to $u$, traverses $\sigma_{u,v}$ (that is, it follows $\pi_{u,v}^{(1)}$, then $\pi_{u,v}$, and then $\pi_{u,v}^{(2)}$), from $u$ to $v$, and then follows $\Gamma$ from $v$ to $y$. Writing $T'$ for the passage time in the weights $(t_e^{u,v})$, we have
\[
T'(\tau) - T(\Gamma) \leq T'(\tau') - T(\Gamma).
\]
The path $\Gamma$ consists of three consecutive subpaths: $\Gamma_1$, which connects $x$ to $u$, $\Gamma_2$, which connects $u$ to $v$, and $\Gamma_3$, which connects $v$ to $y$. Hence $T(\Gamma) = T(\Gamma_1) + T(\Gamma_2) + T(\Gamma_3)$ and similarly $T'(\tau') = T'(\Gamma_1) + T'(\sigma_{u,v}) + T'(\Gamma_3)$. Because $T(\Gamma_2) = T(u,v)$, we obtain from the above display that
\begin{equation}\label{eq: second_path_inequality}
T'(\tau) - T(\Gamma) \leq [ T'(\Gamma_1) - T(\Gamma_1)] + [T'(\sigma_{u,v}) - T(u,v)] + [T'(\Gamma_3)-T(\Gamma_3)].
\end{equation}
Beginning with $\Gamma_1$, we observe that it only contains one vertex in $R(n)$ (the vertex $u$), and therefore has only one edge in $\overline{R}(n)$, the edge that connects $u$ to $R(n)^c$. Similarly, $\Gamma_3$ has only one edge in $\overline{R}(n)$, the edge that connects $v$ to $R(n)^c$. The weights $(t_e)$ and $(t_e^{u,v})$ only differ in $\overline{R}(n)$, and so using occurrence of $A_{u,v}^\ast$, we obtain
\[
T'(\Gamma_1) \leq T(\Gamma_1) + 2M \text{ and } T'(\Gamma_3) \leq T(\Gamma_3) + 2M.
\]
Putting this in \eqref{eq: second_path_inequality} produces
\[
T'(\tau) - T(\Gamma) \leq 4M + T'(\sigma_{u,v}) - T(u,v).
\]
Again using occurrence of $A_{u,v}^\ast$, we can estimate
\[
T'(\sigma_{u,v}) \leq r' \# \pi_{u,v}^{(1)} + \left( r + \frac{\delta}{2}\right)\# \pi_{u,v} + r' \# \pi_{u,v}^{(2)}.
\]
As $\#\pi_{u,v}^{(i)} \leq d-1$ for $i=1,2$, we can combine the previous two displays to get
\begin{equation}\label{eq: crocodile}
T'(\tau) - T(\Gamma) \leq 4M + 2r'(d-1) + \left( r + \frac{\delta}{2} \right) \#\pi_{u,v} - T(u,v).
\end{equation}

To estimate \eqref{eq: crocodile}, we consider two cases. If $\|u-v\|_1 \geq n/2$, then $\overline{\pi}_{u,v}$ is oriented, and so $\#\pi_{u,v} = 2 + \#\overline{\pi}_{u,v} = 2 + \|\overline{u} - \overline{v}\|_1 \leq 2 + 2d + \|u-v\|_1$. From Def.~\ref{def: black}, we get $T(u,v) \geq (r+\delta)\|u-v\|_1$, and if we put this in \eqref{eq: crocodile}, we get
\begin{align}
T'(\tau) - T(\Gamma) &\leq 4M + 2r'(d-1) + (2+2d)\left( r + \frac{\delta}{2}\right)  - \frac{\delta}{2} \|u-v\|_1 \nonumber \\
&\leq 4M + 2r'(d-1) + (2+2d) \left( r + \frac{\delta}{2}\right)- \frac{\delta}{4} n < 0. \label{eq: case_1}
\end{align}
The final inequality follows from our assumption on $n$ in the statement of the proposition, assuming that $M$ is sufficiently large.

The second case is that $\|u-v\|_1 \leq n/2$. Here, we use that $T(u,v) = T(\Gamma_2)$, and because $\Gamma$ crosses $R(n)$ the path $\Gamma_2$ must also cross $R(n)$. So let $a,b$ be two endpoints of a subpath of $\Gamma_2$ that is a crossing of $R(n)$. Because $\|a-b\|_1 \geq n$ and $R(n)$ is black, we obtain $T(u,v) \geq T(a,b) \geq (r+\delta) n$, and if we put this in \eqref{eq: crocodile} and use $\#\pi_{u,v} \leq n+2$, we get 
\begin{align*}
T'(\tau) - T(\Gamma) &\leq 4M + 2r'(d-1) + \left( r+ \frac{\delta}{2}\right) (n+2) - (r+\delta) n \\
&= 4M+2r'(d-1) +2\left( r + \frac{\delta}{2}\right) - \frac{\delta}{2} n.
\end{align*}
Our assumption on $n$ in the statement of the proposition implies that the above quantity is negative if $M$ is sufficiently large. Taking this along with \eqref{eq: case_1}, we conclude that
\[
T'(\tau) - T(\Gamma) < 0 \text{ if } M \text{ is large enough},
\]
and therefore $T'(x,y) < T(x,y)$. Because the weights $(t_e)$ and $(t_e^{u,v})$ only differ on the set $\cup_{i=1}^4 E_i(u,v)$, the geodesic $\tau$ must contain an edge from this union. But it begins at $x$, which is not in $R(n)$, so it must take an edge of $E_3(u,v) \cup E_4(u,v)$ before taking one from $E_1(u,v) \cup E_2(u,v)$. This shows \eqref{eq: take_large_edge}.

To complete the proof of \eqref{eq: main_containment}, we must show that in $(t_e^{u,v})$, the rung $\pi_{u,v}$ is the only rung with at least $n/2$ many edges, all whose edges $e$ have $t_e^{u,v} \leq r + \delta/2$. Item 1 of $A_{u,v}^\ast$ directly implies that $t_e^{u,v} \leq r+\delta/2$ for all edges $e$ of $\pi_{u,v}$, so $\pi_{u,v}$ is one such rung. Also, by definition of $A_{u,v}^\ast$, all edges $e$ that are incident to $\pi_{u,v}$ but not in $\pi_{u,v}$ have $t_e^{u,v} \geq r+\delta$ (such edges are in some $E_i(u,v)$ for $i=2, 3, 4$). So if $\pi'$ is any rung with at least $n/2$ many edges, all whose edges $e$ have $t_e^{u,v} \leq r+\delta/2$, it cannot contain an edge incident to $\pi_{u,v}$. If $\pi'$ starts at $u^\ast$ or $v^\ast$, $\pi'$ must therefore equal $\pi_{u,v}$ or its reversal. If, on the other hand, $\pi'$ did not start at $u^\ast$ or $v^\ast$, then it would have to remain disjoint from $\pi_{u,v}$, else it would contain an edge incident to $\pi_{u,v}$. By definition of $A_{u,v}^\ast$, $\pi'$ could therefore contain no edges of $\cup_{i=1}^4 E_i(u,v)$, and so we would have $t_e = t_e^{u,v}$ for all edges $e$ of $\pi'$. However, as $R(n)$ is black in the configuration $(t_e)$, $R(n)$ does not have any rungs with at least $\Cr{c: black_constant} \log n$ many edges, all whose edges $e$ satisfy $t_e \leq r+\delta/2$. Therefore no such $\pi'$ can exist and this proves \eqref{eq: main_containment}.

We now return to \eqref{eq: my_decomposition} and use independence of the weights $(t_e)$ and $(t_e^\ast)$ along with the estimate \eqref{eq: A_bound} for
\begin{align*}
&\mathbb{P}(\gamma_{x,y} \text{ crosses } R(n) \text{ and } R(n) \text{ is black}) \\
&= \sum_{u,v \in \partial_i R(n)} \mathbb{P}(A_{u,v}^\ast)^{-1} \mathbb{P}(u_0=u,v_0=v, R(n) \text{ is black}, A_{u,v}^\ast) \\
&\leq e^{\Cr{c: A_bound_constant} n} \mathbb{P}(t_e \in [M,2M])^{-2d} \mathbb{P}(t_e \geq M)^{-\Cr{c: A_bound_constant}n} \\
&\times \sum_{u,v \in \partial_i R(n)} \mathbb{P}(u_0=u,v_0=v, R(n) \text{ is black}, A_{u,v}^\ast).
\end{align*}
Next place \eqref{eq: main_containment} into the above display and use the fact that the weights $(t_e)$ and $(t_e^{u,v})$ are identically distributed to obtain
\begin{align}
&e^{-\Cr{c: A_bound_constant} n} \mathbb{P}(t_e \in [M,2M])^{2d} \mathbb{P}(t_e \geq M)^{\Cr{c: A_bound_constant}n} \mathbb{P}(\gamma_{x,y} \text{ crosses } R(n) \text{ and } R(n) \text{ is black}) \nonumber \\
\leq~& \sum_{u,v \in \partial_i R(n)} \mathbb{P}\left( 
\begin{array}{c}
\text{each geodesic from }x \text{ to } y \text{ contains an edge} \\
e \in \overline{R}(n) \text{ with }t_e \geq M \text{ and } \pi_{u,v} \text{ is the only} \\
\text{rung with at least } \frac{n}{2} \text{ many edges,} \\
\text{all whose edges }e \text{ have } t_e \leq r+\frac{\delta}{2}
\end{array}
\right) \nonumber\\
=~& \sum_{u,v \in \partial_i R(n)} \left( \sum_\pi \mathbf{1}_{\{\pi_{u,v} = \pi\}} \right) \mathbb{P}\left( 
\begin{array}{c}
\text{each geodesic from }x \text{ to } y \text{ contains an edge} \\
e \in \overline{R}(n) \text{ with }t_e \geq M \text{ and } \pi_{u,v} \text{ is the only} \\
\text{rung with at least } \frac{n}{2} \text{ many edges,} \\
\text{all whose edges }e \text{ have } t_e \leq r+\frac{\delta}{2}
\end{array}
\right) \nonumber \\
=~&\sum_{\pi} \bigg[ \#\{(u,v) \in \partial_iR(n) \times \partial_i R(n) : \pi_{u,v} = \pi\} \nonumber \\
\times~&\mathbb{P}\left( 
\begin{array}{c}
\text{each geodesic from }x \text{ to } y \text{ contains an edge} \\
e \in \overline{R}(n) \text{ with }t_e \geq M \text{ and } \pi \text{ is the only} \\
\text{rung with at least } \frac{n}{2} \text{ many edges,} \\
\text{all whose edges }e \text{ have } t_e \leq r+\frac{\delta}{2}
\end{array}
\right) \bigg]. \label{eq: big_equation}
\end{align}
The sum is over all rungs $\pi$ of $R(n)$ with at least $n/2$ many edges. The number of pairs $(u,v)$ for which $\pi_{u,v} = \pi$ can be estimated as follows. Letting $a$ and $b$ be the endpoints of $\pi$, if $(u,v)$ is such that $\pi_{u,v} = \pi$, then the sets $\{u^\ast,v^\ast\}$ and $\{a,b\}$ must be equal. Therefore
\begin{align*}
\#\{(u,v) : \pi_{u,v} = \pi\} &\leq \#\{(u,v) : \{u^\ast, v^\ast\} = \{a,b\}\} \\
&\leq \#\{u : u^\ast \in \{a,b\}\} \times \#\{v : v^\ast \in \{a,b\}\} \\
&\leq \left( \#\{u : \|u-a\|_\infty \leq d-1\} + \#\{u : \|u-b\|_\infty \leq d-1\} \right)^2 \\
&\leq 4(2d-1)^{2d}.
\end{align*}
This implies the bound
\begin{equation}\label{eq: what_we_need_for_chris}
\eqref{eq: big_equation} \leq 4(2d-1)^{2d} \sum_\pi \mathbb{P}\left( 
\begin{array}{c}
\text{each geodesic from }x \text{ to } y \text{ contains an edge} \\
e \in \overline{R}(n) \text{ with }t_e \geq M \text{ and } \pi \text{ is the only} \\
\text{rung with at least } \frac{n}{2} \text{ many edges,} \\
\text{all whose edges }e \text{ have } t_e \leq r+\frac{\delta}{2}
\end{array}
\right).
\end{equation}
As $\pi$ varies, the events in the probability are disjoint, so
\[
\eqref{eq: big_equation}\leq 4(2d-1)^{2d} \mathbb{P}(\text{each geodesic from }x \text{ to } y \text{ contains an edge } e \in \overline{R}(n) \text{ with } t_e \geq M).
\]
Combining this with the display leading to \eqref{eq: big_equation} and taking $M$ so large that $(4(2d-1))^{-2d} e^{-\Cr{c: A_bound_constant}} \geq \mathbb{P}(t_e \geq M)$ completes the proof of the proposition.
\end{proof}

\medskip
\noindent
{\bf Step 3: Controlling dependence using good rectangles.}
Now that we have proved Prop.~\ref{prop: modification}, we must estimate the number of black rectangles that geodesics cross. The conditions in Def.~\ref{def: black} are not finitely dependent (specifically $T(a,b)$ in item 2 depends on weights arbitrarily far outside of $R(n)$), so we need to ensure geodesics that begin and end in $R(n)$ do not travel too far from $R(n)$.
\begin{df}
Let $\Cl[lgc]{c: good_constant}\geq 3$ and $n \geq 1$. We say that the rectangle $R(n)$ is \underline{good} if for each $w_1,w_2 \in R(n)$, we have
\[
T_{B(\Cr{c: good_constant}n)}(w_1,w_2) < T(w_1, \partial B(\Cr{c: good_constant}n)) + T(w_2, \partial B(\Cr{c: good_constant}n)),
\]
where $T_{B(\Cr{c: good_constant}n)}(w_1,w_2)$ is the minimal passage time from $w_1$ to $w_2$ among all paths from $w_1$ to $w_2$ that remain in $B(\Cr{c: good_constant}n)$, $B(\Cr{c: good_constant}n) = [-\Cr{c: good_constant} n, \Cr{c: good_constant}n]^d$, and
\[
\partial B(\Cr{c: good_constant}n) = \{x \in \mathbb{Z}^d \setminus B(\Cr{c: good_constant}n) : \exists~y \in B(\Cr{c: good_constant}n) \text{ with } |x-y| = 1\}.
\]
\end{df}

Note that if $R(n)$ is good, then any geodesic between any two vertices $w_1,w_2 \in R(n)$ must remain in $B(\Cr{c: good_constant}n)$. Therefore
\begin{equation}\label{eq: finite_dependence}
\{R(n) \text{ is black and good}\} \text{ depends only on edge-weights in } B(\Cr{c: good_constant} n).
\end{equation}
We now show that if $n$ is large, then $R(n)$ is good with high probability.
\begin{lem}\label{lem: another_good_lemma}
There exists $\Cr{c: good_constant}>0$ such that
\[
\mathbb{P}(R(n) \text{ is good}) \to 1 \text{ as } n \to \infty.
\]
\end{lem}
\begin{proof}
The proof of this lemma is similar to that of Lem.~\ref{lem: good_lemma}, and once again uses percolation methods from the proof of Prop.~\ref{prop: improved_tail_bound}. Let $M_0$ be such that $F(M_0) \geq p_0$ (from Lem.~\ref{lem: gamma_lemma}). We call an edge $e$ \underline{open} if $t_e \leq M_0$; otherwise, it is \underline{closed}. We again let $I$ be the unique infinite open cluster but now define 
\[
A_n = \{\text{every path from } \partial B(3n) \text{ to } \partial B(4n) \text{ intersects }I\}.
\]
By Lem.~\ref{lem: gamma_lemma}, we have
\begin{align*}
\mathbb{P}(A_n^c) &\leq \sum_{w \in \partial B(3n)} \mathbb{P}(\text{some path from } w \text{ to } w + \partial B(n) \text{ does not intersect }I) \\
&\leq (6n+1)^d e^{-\Cr{c: book_chemical_constant}n} \to 0 \text{ as } n \to \infty.
\end{align*}

Similarly to \eqref{eq: prime_bound}, if we define
\[
A_n' = \{\forall z_1,z_2 \in I \cap B(5n), \text{ we have } d_I(z_1,z_2) \leq \Cl[lgc]{c: new_antal_pisztora}n\}
\]
for some constant $\Cr{c: new_antal_pisztora}>0$ (here $d_I(z_1,z_2)$, as before, is the smallest number of edges in any open path from $z_1$ to $z_2$), then there exists $\Cr{c: new_antal_pisztora}$ such that
\[
\mathbb{P}(A_n') \to 1 \text{ as } n \to \infty.
\]

Last, we define
\[
B_n = \{T(\partial B(4n), \partial B(\Cr{c: good_constant}n)) \geq a(\Cr{c: good_constant} - 4)n\},
\]
where $a$ is from the statement of Lem.~\ref{lem: kesten_lemma}. Using Lem.~\ref{lem: kesten_lemma} and a union bound over vertices in $\partial B(4n)$ gives for any $\Cr{c: good_constant}\geq 5$ that
\[
\mathbb{P}(B_n^c) = \mathbb{P}(T(\partial B(4n), \partial B(\Cr{c: good_constant}n)) < a (\Cr{c: good_constant}-4)n) \leq (4n+1)^d e^{-\Cr{c: kesten_lemma_constant} (\Cr{c: good_constant}-4)n} \to 0 \text{ as } n \to \infty.
\]

We now choose $\Cr{c: good_constant} \geq 5$ to be sufficiently large. Specifically, it should satisfy
\begin{equation}\label{eq: constant_condition}
\Cr{c: good_constant} \geq \Cr{c: new_antal_pisztora} + 4 \text{ and } 2a(\Cr{c: good_constant} - 4) > \Cr{c: new_antal_pisztora}M_0.
\end{equation}

See Fig.~\ref{fig: fig_4} for an illustration of the following argument. Choose any outcome in $A_n \cap A_n' \cap B_n$. 
For $w_1,w_2 \in R(n)$, let $\gamma_1$ and $\gamma_2$ be optimal paths for $T(w_1, \partial B(\Cr{c: good_constant}n))$ and $T(w_2, \partial B(\Cr{c: good_constant}n))$, respectively. By occurrence of $A_n$, we can pick a vertex $c \in \gamma_1 \cap I$ that is on the initial segment of $\gamma_1$ from $w_1$ to $\partial B(4n)$, and we can pick a vertex $d \in \gamma_2 \cap  I$ that is on the initial segment of $\gamma_2$ from $w_2$ to $\partial B(4n)$.
We can apply the condition in $A_n'$ to the pair $(c,d)$ to obtain 
\begin{equation}\label{eq: pi_length}
d_I(c,d) \leq \Cr{c: new_antal_pisztora}n.
\end{equation}
Therefore we can pick a path $\gamma_3$ from $c$ to $d$ that remains in $I$ and has at most $\Cr{c: new_antal_pisztora}n$ many edges.

\begin{figure}[p]
  \centering
  \includegraphics[width=.9\linewidth,trim={0cm 5cm 0cm 3.5cm}, clip]{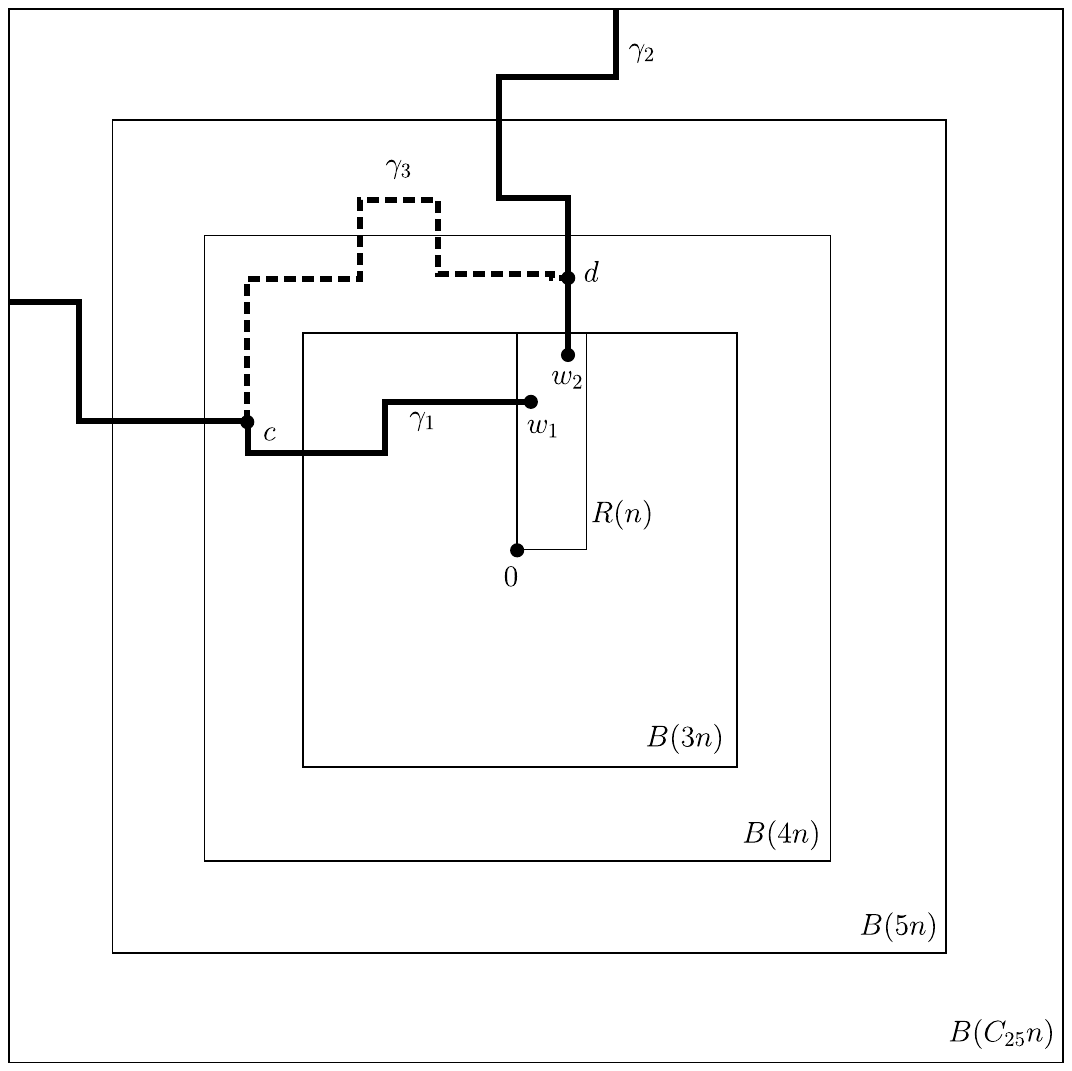}
  \caption{Illustration of the argument for Lem.~\ref{lem: another_good_lemma}. The paths $\gamma_1$ and $\gamma_2$ are optimal paths for $T(w_1,\partial B(\Cr{c: good_constant}n))$ and $T(w_2,\partial B(\Cr{c: good_constant}n))$. The path $\gamma_3$ is contained in the infinite component $I$ and connects points $c$ on $\gamma_1$ and $d$ on $\gamma_2$. Because the union of $\gamma_3$ with the initial segments of $\gamma_1$ and $\gamma_2$ comprises a path from $w_1$ to $w_2$ with smaller passage time than the sum of $T(\gamma_1)$ and $T(\gamma_2)$, no geodesic from $w_1$ to $w_2$ may exit $B(\Cr{c: good_constant}n)$.}
  \label{fig: fig_4}
\end{figure}

By occurrence of $B_n$, we have 
\[
T(c, \partial B(\Cr{c: good_constant}n)) \geq a(\Cr{c: good_constant}-4)n \text{ and } T(d, \partial B(\Cr{c: good_constant}n)) \geq a(\Cr{c: good_constant}-4)n,
\]
and so
\begin{align*}
&T(w_1, \partial B(\Cr{c: good_constant}n)) + T(w_2, \partial B(\Cr{c: good_constant}n)) \\
=~& T(w_1,c) + T(w_2,d) + T(c, \partial B(\Cr{c: good_constant}n)) + T(d, \partial B(\Cr{c: good_constant}n)) \\
\geq~& T(w_1,c) + T(w_2,d) + 2a(\Cr{c: good_constant}-4)n.
\end{align*}
Let $\pi$ be the path obtained by starting at $w_1$, following $\gamma_1$ to $c$, following $\gamma_3$ to $d$, and then following $\gamma_2$ to $w_2$. Then by \eqref{eq: pi_length} and the choice of $c$ and $d$, $\pi$ remains in the box $B((\Cr{c: new_antal_pisztora} + 4) n) \subset B(\Cr{c: good_constant}n)$ (the inclusion holds by the first inequality of \eqref{eq: constant_condition}) and so
\[
T_{B(\Cr{c: good_constant}n)}(w_1,w_2) \leq T(\pi) \leq T(w_1,c) + T(w_2,d) + \Cr{c: new_antal_pisztora}M_0 n.
\]
Combining the previous two displays and using the second inequality of \eqref{eq: constant_condition} gives 
\[
T_{B(\Cr{c: good_constant}n)}(w_1,w_2) < T(w_1, \partial B(\Cr{c: good_constant}n)) + T(w_2, \partial B(\Cr{c: good_constant}n))
\]
and therefore proves that
\[
A_n \cap A_n' \cap B_n \subset \{R(n) \text{ is good}\}.
\]
Because $\mathbb{P}(A_n \cap A_n' \cap B_n) \to 1$, this completes the proof.
\end{proof}

\medskip
\noindent
{\bf Step 4: Finishing the proof.}
From now on, we fix $\Cr{c: good_constant}$ such that the conclusion of Lem.~\ref{lem: another_good_lemma} holds. In this last part of the proof, we put together the above pieces with a block argument modeled off of the one in \cite{BK93}. Let $x \in \mathbb{Z}^d$; we will work with the geodesic $\gamma_{0,x}$ and estimate the number of good black rectangles that it crosses. We will use the terms ``good'' and ``black'' not only when we refer to the rectangle $R(n)$, but also in an analogous way when we refer to any integer translation or rotation of $R(n)$. For $n \geq 1$, define the family of blocks
\[
\{S(n,z) : z \in \mathbb{Z}^d\}, \text{ where } S(n,z) = nz + [0,n]^d
\]
and the family of larger blocks
\[
\left\{\overline{S}(n,z) : z \in \mathbb{Z}^d\right\}, \text{ where } \overline{S}(n,z) = nz + [-n,2n]^d.
\]
The annulus $(\overline{S}(n,z) \setminus S(n,z)) \cup \partial_i S(n,z)$ can be covered by $2d$ overlapping copies of the rectangle $R(n)$ in the natural way, by using integer translations and rotations of $R(n)$. We enumerate these rectangles in any deterministic fashion as $R_1(n,z), \dots, R_{2d}(n,z)$.
\begin{df}\label{def: black_cube}
We say that the cube $S(n,z)$ is \underline{black} (respectively \underline{good}) if all of the rectangles $R_1(n,z), \dots, R_{2d}(n,z)$ are black (respectively good).
\end{df}
Because of \eqref{eq: finite_dependence}, the variables $(\mathbf{1}_{\{S(n,z) \text{ is black and good}\}} : z \in \mathbb{Z}^d)$ form a finitely dependent site percolation process. Furthermore, by Lemmas~\ref{lem: black} and \ref{lem: another_good_lemma}, for each $z$, we have $\mathbb{P}(S(n,z) \text{ is black and good}) \to 1$ as $n \to \infty$. We can then apply a standard Peierls argument (see \cite[Lem.~5.2]{BK93} and \cite[proof of (3.12)]{GK84}) to deduce that there exists $\Cl[smc]{c: vdb_block_constant} > 0$ such that for all $x \in \mathbb{Z}^d$ and all sufficiently large $n$, we have
\[
\mathbb{P}\left(\exists \text{ a path from } 0 \text{ to } x \text{ that visits at most } \frac{\Cr{c: vdb_block_constant}}{n} \|x\|_1 \text{ distinct good black cubes}\right) \leq e^{-\Cr{c: vdb_block_constant} \frac{\|x\|_1}{n}}.
\]
We henceforth fix $n = \left\lceil 64 M / \delta \right\rceil$ and take $M$ sufficiently large so that both the result of Prop.~\ref{prop: modification} and the previous display hold. By applying the display, we have for some $\Cl[smc]{c: good_black_number}, \Cl[lgc]{c: my_new_constant_now}>0$, all sufficiently large $M$, and all $x \in \mathbb{Z}^d$ with $\|x\|_1 \geq \Cr{c: my_new_constant_now} M$,
\[
\mathbb{E}\#\{\text{distinct black cubes that } \gamma_{0,x} \text{ visits}\} \geq \frac{\Cr{c: good_black_number}}{M} \|x\|_1.
\]
If $\gamma_{0,x}$ visits a cube $S(n,z)$ but $0$ and $x$ are not in $\overline{S}(n,z)$, then $\gamma_{0,x}$ must contain a crossing of at least one of the $2d$ rectangles $R_i(n,z)$ surrounding $S(n,z)$. As there is at most a constant number of cubes $\overline{S}(n,z)$ that contain $0$ or $x$, we use these comments and the above display to find $\Cl[lgc]{c: dimension_constant}>0$ such that for all $x \in \mathbb{Z}^d$ and all sufficiently large $M$, we have
\begin{align*}
\frac{\Cr{c: good_black_number}}{M}\|x\|_1 - \Cr{c: dimension_constant} &\leq \sum_{z : 0,x \notin \overline{S}(n,z)} \mathbb{P}(\gamma_{0,x} \text{ visits } S(n,z) \text{ and } S(n,z) \text{ is black}) \\
&\leq \sum_{z : 0,x \notin \overline{S}(n,z)} \sum_{i=1}^{2d} \mathbb{P}(\gamma_{0,x} \text{ crosses } R_i(n,z) \text{ and } R_i(n,z) \text{ is black}).
\end{align*}
Now we multiply both sides of this inequality by the two prefactor probabilities in Prop.~\ref{prop: modification}, and then apply this proposition to each rectangle to obtain
\begin{align*}
&\mathbb{P}(t_e \in [M,2M])^{2d} \mathbb{P}(t_e \geq M)^{\Cr{c: modification_constant} n}\left( \frac{\Cr{c: good_black_number}}{M}\|x\|_1 - \Cr{c: dimension_constant}\right) \nonumber \\
\leq~& \sum_{z : 0,x \notin \overline{S}(n,z)} \sum_{i=1}^{2d} \mathbb{P}\left(
\begin{array}{c}
\text{each geodesic from } 0 \text{ to } x \text{ contains} \\
\text{an edge }e \in \overline{R}_i(n,z) \text{ with } t_e \geq M
\end{array}
\right) \\
\leq~& \sum_{z : 0,x \notin \overline{S}(n,z)} \sum_{i=1}^{2d} \mathbb{E}\inf_\gamma \#\{e \in \overline{R}_i(n,z) : e \in \gamma \text{ and } t_e \geq M\} \\
\leq~& \mathbb{E}\inf_\gamma \left( \sum_{z : 0,x \notin \overline{S}(n,z)} \sum_{i=1}^{2d} \sum_e \mathbf{1}_{\{e \in \overline{R}_i(n,z)\}} \mathbf{1}_{\{e \in \gamma \text{ and } t_e \geq M\}}\right) \\
\leq~& (2d)5^d \mathbb{E} \inf_\gamma \sum_e \mathbf{1}_{\{e \in \gamma \text{ and } t_e \geq M\}}.
\end{align*}
In the last three lines, the infimum is over all geodesics $\gamma$ from $0$ to $x$. In the fourth inequality, we have used that the number of sets $\overline{R}_i(n,z)$ that can contain a given edge is at most $(2d)5^d$ (at most $5^d$ choices for $z$ and $2d$ choices for $i$). Because the final expectation equals $\mathbb{E}N_{\mathrm{min}}(0,x,[M,\infty))$, we can complete the proof by choosing the value of $n = \left\lceil 64M/\delta \right\rceil$.

\section{Bounds for general sets $A$}\label{sec: general_A}

\subsection{Proof of Prop.~\ref{prop: DHS_improvement}}
Recalling that $\textsc{GEO}(0,x)$ is the edge self-avoiding geodesic from $0$ to $x$ with the most edges, we decompose
\begin{equation}\label{eq: first_breakdown_4_prime}
\mathbb{E}N_{\mathrm{max}}(0,x,A) = \sum_{k=0}^\infty \mathbb{E}N_{\mathrm{max}}(0,x,A) \mathbf{1}_{\{\#\textsc{GEO}(0,x) \in [k|x|,(k+1)|x|)\}}.
\end{equation}

On the event $\{\#\textsc{GEO}(0,x) < (k+1)|x|\}$, we have $N_{\mathrm{max}}(0,x,A) \leq N_{(k+1)|x|}$ almost surely, where $N_n$ is defined as follows (and in Appendix~\ref{sec: appendix}). Set $X_e = 1$ if $t_e \in A$ and $X_e=0$ otherwise. Then $N_n$ is the maximal value of $\sum_{e \in S} X_e$ over all subsets $S$ of $\mathcal{E}^d$ of cardinality $n$ which contain an edge incident to the origin and are connected. Combining this with the Cauchy-Schwarz inequality gives
\begin{equation}\label{eq: first_breakdown_4_prime_again}
\mathbb{E}N_{\mathrm{max}}(0,x,A) \mathbf{1}_{\{\#\textsc{GEO}(0,x) \in [k|x|,(k+1)|x|)\}} \leq \sqrt{\mathbb{E}N_{(k+1)|x|}^2} \sqrt{ \mathbb{P}(\#\textsc{GEO}(0,x) \geq k|x|)}.
\end{equation}
By \eqref{eq: lattice_animals_second_moment}, we have
\begin{equation}\label{eq: LA_to_put_in}
\mathbb{E}N^2_{(k+1)|x|} \leq \left( \Cr{c: 2024_constant_1} \mathbb{P}(t_e \in A)^{\frac{1}{d}} (k+1)|x|\right)^2
\end{equation}
Also, Prop.~\ref{prop: improved_tail_bound} implies the much weaker estimate
\[
\mathbb{P}(\#\textsc{GEO}(0,x) \geq k|x|) \leq \frac{\Cl[lgc]{c: chaz_supreme}}{(k+1)^5}
\]
for some $\Cr{c: chaz_supreme}$ independent of $x$ and $k$. We put this and \eqref{eq: LA_to_put_in} into the right side of \eqref{eq: first_breakdown_4_prime_again} and then back into the right side of \eqref{eq: first_breakdown_4_prime} to obtain, for some $\Cr{c: my_last_constant}>0$,
\[
\mathbb{E}N_{\mathrm{max}}(0,x,A) \leq \sum_{k=0}^\infty \Cr{c: 2024_constant_1} \mathbb{P}(t_e \in A)^{\frac{1}{d}} (k+1)|x| \sqrt{\frac{\Cr{c: chaz_supreme}}{(k+1)^5}} \leq \Cl[lgc]{c: my_last_constant} \mathbb{P}(t_e \in A)^{\frac{1}{d}}|x|.
\]

\subsection{Proof of Prop.~\ref{prop: separation_prop}}
Choose $a \in (r,b)$ such that $F(a)(F(b)-F(a))>0$; this is possible because $F(b)>F(r)$ and $r$ is the infimum of the support of $F$. We will use the variable $D_e$ which was defined in \eqref{eq: D_e_def} as $D_e = \sup\{s \geq 0 : e \in \overline{\textsc{GEO}}(0,x) \text{ in } (s,t_e^c)\}$, where $t_e^c$ represents the variables $(t_f)_{f \neq e}$ and $s$ is the value of the weight at edge $e$. (Here, if the set is empty, we used the convention $D_e = -\infty$.) In addition, we recall that $t_e \leq D_e$ if and only if $e \in \overline{\textsc{GEO}}(0,x)$.

We then estimate
\[
\mathbb{E}\overline{N}(0,x,(c,d]) = \sum_{e \in \mathcal{E}^d} \mathbb{P}(t_e \leq D_e, t_e \in (c,d]).
\]
By independence of $D_e$ and $t_e$, the above equals
\begin{align}
&\sum_{e \in \mathcal{E}^d} \mathbb{E}(F(\min\{d,D_e\}) - F(\min\{c,D_e\})) \nonumber \\
=~& \sum_{e \in \mathcal{E}^d} \mathbb{E}(F(\min\{d,D_e\}) - F(\min\{c,D_e\}))\mathbf{1}_{\{D_e > c\}} \label{eq: that_one_line}\\
\leq~& (F(d)-F(c)) \sum_{e \in \mathcal{E}^d} \mathbb{P}(D_e > c) \nonumber \\
=~& \frac{F(d)-F(c)}{F(b)-F(a)} \sum_{e \in \mathcal{E}^d} \mathbb{E}(F(\min\{b,D_e\})-F(\min\{a,D_e\}))\mathbf{1}_{\{D_e > c\}} \nonumber \\
\leq~& \frac{F(d)-F(c)}{F(b)-F(a)} \sum_{e \in \mathcal{E}^d} \mathbb{E}(F(\min\{b,D_e\})-F(\min\{a,D_e\}))\mathbf{1}_{\{D_e > a\}}. \nonumber
\end{align}
For the second equality, we have used that $a \leq b \leq c$. By running the argument that led to \eqref{eq: that_one_line} in reverse, with $a,b$ in place of $c,d$, the sum in the last line equals $\mathbb{E}\overline{N}(0,x,(a,b])$. We obtain
\[
\mathbb{E}\overline{N}(0,x,(c,d]) \leq \frac{F(d)-F(c)}{F(b)-F(a)} \mathbb{E}\overline{N}(0,x,(a,b]).
\]
Applying Thm.~\ref{thm: upper_bound_upper_tail} (with $M=a$ and bounding $q(a,x)$ by the expression of the first item in that theorem) to the right side completes the proof:
\begin{align*}
\mathbb{E}\overline{N}(0,x,(c,d]) &\leq 4 \Cr{c: L1} \frac{F(d)-F(c)}{F(b)-F(a)} \cdot \frac{1-F(a)}{F(a)} |x| \\
&\leq \frac{4 \Cr{c: L1}}{F(a)(F(b)-F(a))} \cdot (F(d)-F(c))|x|.
\end{align*}

\subsection{Proof of Prop.~\ref{prop: FKG_lower_bound}}

To prove Prop.~\ref{prop: FKG_lower_bound}, we again use the variable $D_e$ to estimate
\begin{align*}
\mathbb{E}\underline{N}(0,x,A) &\geq \sum_{e \in \mathcal{E}^d} \mathbb{P}(t_e \leq a \text{ and } t_e < D_e) \\
&= \sum_{e \in \mathcal{E}^d} \lim_{\epsilon \to 0} \mathbb{P}(t_e \leq a \text{ and } t_e \leq D_e - \epsilon).
\end{align*}
Because $D_e$ is a function of $t_e^c$, which is independent of $t_e$, the probability in the sum above equals $\mathbb{E} F(\min\{a,D_e-\epsilon\})$. However, if $F$ is any distribution function, we have for all $c,d \in \mathbb{R}$,
\[
F(\min\{c,d\}) \geq F(c)F(d),
\]
so we obtain the lower bound
\begin{equation} \label{eq: almost_to_nakajima}
\mathbb{E}\underline{N}(0,x,A) \geq F(a) \sum_{e \in \mathcal{E}^d} \mathbb{E}\lim_{\epsilon \to 0} F(D_e - \epsilon) =  F(a) \sum_{e \in \mathcal{E}^d} \mathbb{P}(t_e < D_e). 
\end{equation}

From \eqref{eq: almost_to_nakajima}, we are left to prove that 
\begin{equation}\label{eq: intersection_lower_bound}
\sum_{e \in \mathcal{E}^d} \mathbb{P}(t_e < D_e) \geq \Cr{c: easy_lower_bound} |x|.
\end{equation}
This was essentially done in \cite[Proof of Thm.~3]{N18} with different notation, but the author assumed the moment bound $\mathbb{E}t_e < \infty$ to apply the comparison theorem of van den Berg and Kesten in \cite{BK93}. Here we explain how to use Thm.~\ref{thm: BK_argument} instead to remove this moment condition.

Just as in \cite[Lem.~7]{N18}, we first show that for any $\alpha \in \mathbb{R}$ and any $x \in \mathbb{Z}^d$, we have
\begin{equation}\label{eq: almost_end_of_nakajima}
\sum_{e \in \mathcal{E}^d} \mathbb{P}(t_e < D_e) \geq F(\alpha) \mathbb{E}\overline{N}(0,x,(\alpha,\infty)).
\end{equation}
To do this, let $e_0$ be any edge and let $\hat{t}_{e_0}$ be an independent copy of $t_{e_0}$. Define the configuration $(t_e^\ast)$ by 
\[
t_e^\ast = \begin{cases}
\hat{t}_e & \quad\text{if } e = e_0 \\
t_e & \quad\text{if } e \neq e_0.
\end{cases}
\]
Then, using independence and the fact that the variable $D_{e_0}$ is the same in both configurations $(t_e)$ and $(t_e^\ast)$, we obtain
\begin{align*}
\mathbb{P}(t_{e_0} \leq D_{e_0}, t_{e_0} > \alpha) F(\alpha) &= \mathbb{P}(t_{e_0} \leq D_{e_0}, t_{e_0} > \alpha, \hat{t}_{e_0} \leq \alpha) \\
&\leq \mathbb{P}(\hat{t}_{e_0} < D_{e_0}) \\
&= \mathbb{P}(t_{e_0} < D_{e_0}).
\end{align*}
Summing this inequality over $e_0$ (and replacing $e_0$ by $e$) produces
\[
\sum_{e \in \mathcal{E}^d} \mathbb{P}(t_e < D_e) \geq F(\alpha) \sum_{e \in \mathcal{E}^d} \mathbb{P}(t_e \leq D_e, t_e > \alpha),
\]
and this implies \eqref{eq: almost_end_of_nakajima}.

If the distribution of $t_e$ is bounded, then $\mathbb{E}t_e<\infty$, and, as stated in \cite{N18}, the comparison theorem of van den Berg and Kesten applies to show that for some $\alpha>r$, there exists $\Cl[smc]{c: BK_squared}>0$ such that for all $x \in \mathbb{Z}^d$, we have
\[
\mathbb{E}\overline{N}(0,x,(\alpha,\infty)) \geq \Cr{c: BK_squared}|x|.
\]
Along with \eqref{eq: almost_end_of_nakajima}, this shows \eqref{eq: intersection_lower_bound} and completes the proof. Otherwise, the distribution of $t_e$ is unbounded. Fixing $\alpha > r$ such that $\mathbb{P}(t_e \in [\alpha + 1, 2(\alpha+1)])>0$ and applying Thm.~\ref{thm: BK_argument} with $\alpha = M$ gives that for $x \in \mathbb{Z}^d$ with $|x|$ sufficiently large, we have
\[
\mathbb{E}\overline{N}(0,x,(\alpha,\infty)) \geq \mathbb{E}N_{\mathrm{min}}(0,x,[\alpha+1,\infty)) \geq \Cl[smc]{c: new_nakajima} |x|.
\]
Combined with \eqref{eq: almost_end_of_nakajima}, this shows \eqref{eq: intersection_lower_bound} and completes the proof if $|x|$ is sufficiently large.

To show \eqref{eq: intersection_lower_bound} in the case of small $|x|$ but unbounded distribution for $t_e$, we simply need to show (because of \eqref{eq: almost_end_of_nakajima}) that for each nonzero $x$, there exists $\alpha > r$ such that $\mathbb{E}\overline{N}(0,x,(\alpha,\infty)) > 0$. To this end, let $x \in \mathbb{Z}^d$ be nonzero and pick any (say) oriented path $\gamma_x$ from $0$ to $x$. Choose $\alpha > r$ such that $\mathbb{P}(t_e \in (\alpha,\alpha+1)) > 0$ and let $\beta > (\alpha+1)\|x\|_1$. Define $A_x$ to be the event that (a) each edge $e$ in $\gamma_x$ has $t_e \in (\alpha,\alpha+1)$ and (b) each edge $e$ that is adjacent to, but not in, $\gamma_x$ has $t_e > \beta$. Then $\mathbb{P}(A_x) > 0$ but we claim that for any outcome in $A_x$, $\gamma_x$ is a geodesic from $0$ to $x$ all whose edges have weight greater than $\alpha$. Indeed, for such an outcome, if $\Gamma$ is any path from $0$ to $x$ that is not equal to $\gamma_x$, then  it must contain an edge adjacent to, but not in, $\gamma_x$ and so
\[
T(\Gamma) > \beta > (\alpha+1)\|x\|_1 > T(\gamma_x).
\]
Thus $\mathbb{E}\overline{N}(0,x,(\alpha,\infty)) \geq \mathbb{P}(A_x)>0$ and this completes the proof.

\appendix

\section{Lattice animal bound}\label{sec: appendix}

\begin{df}\label{def: k_dependent}
Let $p \in (0,1)$ and let $(X_e)_{e \in \mathcal{E}^d}$ be a family of random variables satisfying
\[
\mathbb{P}(X_e = 1) = p = 1-\mathbb{P}(X_e=0) \text{ for all }e.
\]
We say that $(X_e)$ is \underline{$k$-dependent} if whenever $\mathcal{E}_1,\mathcal{E}_2 \subset \mathcal{E}^d$ satisfy $\min\{d(e,f) : e \in \mathcal{E}_1,f \in \mathcal{E}_2\} > k$, the collections $(X_e)_{e \in \mathcal{E}_1}$ and $(X_e)_{e \in \mathcal{E}_2}$ are independent. Here, $d(e,f)$ is the minimal $\ell^\infty$-distance between any endpoint of $e$ and any endpoint of $f$.
\end{df}

For $n \geq 1$, let $\Xi_n$ be the collection of subsets of $\mathcal{E}^d$ of cardinality $n$ which contain an edge incident to the origin and which are connected (a set $B$ of edges is connected if for each $e,f \in B$, there is a path with initial edge $e$ and final edge $f$ all whose edges lie in $B$). A \underline{lattice animal} is a finite connected subgraph of $\mathbb{Z}^d$ containing the origin, and the greedy lattice animal problem with weights $(X_e)$ seeks to estimate the maximal weight of any lattice animal of size $n$. Accordingly, we define
\[
N_n = \max_{S \in \Xi_n} \left( \sum_{e \in S} X_e\right).
\]
The main result of this appendix is the following extension of \cite[Lem.~6.8]{DHS15}, which itself extends  \cite{L97} and \cite[Prop.~2.2]{M02}.
\begin{prop}\label{prop: lattice_animals}
There exists $\Cr{c: lattice_animal_constant}>0$ such that for all $p \in (0,1)$, all $k,n \geq 1$, and all families $(X_e)_{e \in \mathcal{E}^d}$ of $k$-dependent Bernoulli$(p)$ random variables,
\[
\mathbb{E}N_n \leq \Cl[lgc]{c: lattice_animal_constant} k^d p^{\frac{1}{d}}n.
\]
\end{prop}
\begin{proof}
If $n \leq p^{-1/d}$, then we have
\[
\mathbb{E}\left(\frac{N_n}{p^{\frac{1}{d}}n}\right) \leq \frac{1}{p^{\frac{1}{d}}n} \sum_{e \in [-n,n]^d} \mathbb{E}X_e \leq \frac{2d(2n+1)^dp}{p^{\frac{1}{d}}n} \leq 2d(3^d) \left( p^{\frac{1}{d}}n\right)^{d-1} \leq 3^{d+1}d.
\]

From now, we suppose that $n > p^{-\frac{1}{d}}$. To estimate $N_n$, we will cover any element of $\Xi_n$ by cubes of sidelength $p^{-1/d}$. We use \cite[Lem.~1]{CGGK93}, which states that for any connected $A \subset \mathbb{Z}^d$ containing the origin and having cardinality $n$, and each integer $\ell \in [1,n]$, one can find vertices $x_0, \dots, x_r \in \mathbb{Z}^d$ with $r = \lfloor 2n/\ell\rfloor$ such that
\begin{enumerate}
\item[(1)] $x_0=0$ and $\|x_{i+1}-x_i\|_\infty \leq 1$ for $i=0, \dots, r-1$, and
\item[(2)] $A \subset \cup_{i=0}^r \left(\ell x_i + [-\ell,\ell]^d\right)$.
\end{enumerate}
For any $S \in \Xi_n$, the set of endpoints of edges of $S$ is connected, contains the origin, and has cardinality at most $n+1$. Using this result with $\ell = \lceil p^{-1/d}\rceil$ and a union bound, we obtain
\begin{equation}\label{eq: overall_N_n}
\mathbb{P}\left( N_n \geq snp^{\frac{1}{d}}\right) \leq \sum_{x_0, \dots, x_r} \mathbb{P}\left( \sum_{e \in \cup_{i=0}^r (\ell x_i + [-\ell,\ell]^d)} X_e \geq snp^{\frac{1}{d}}\right).
\end{equation}
The subscript of the inner sum means that both endpoints of $e$ are in the union. Also, the outer sum is over all choices of  sequences $x_0, \dots, x_r$ satisfying item (1) above. The number of terms in such a sequence is 
\[
r+1 = 1 + \left\lfloor \frac{2(n+1)}{\ell}\right\rfloor \leq 5np^{\frac{1}{d}}
\]
and the number of elements in the innermost sum is
\[
\#\left\{e \in \cup_{i=0}^r (\ell x_i + [-\ell,\ell]^d)\right\} \leq(r+1)(2d)(2\ell+1)^d \leq \Cl[lgc]{c: appendix_1} np^{\frac{1}{d}-1}
\]
for some $\Cr{c: appendix_1}>0$.

For any sequence $\mathbf{x} = (x_0, \dots, x_r)$, define
\[
U_{\mathbf{x}} = \sum_{e \in \cup_{i=0}^r (\ell x_i + [-\ell,\ell]^d)} X_e,
\]
so that \eqref{eq: overall_N_n} becomes
\[
\mathbb{P}\left( N_n \geq snp^{\frac{1}{d}}\right) \leq \sum_{\mathbf{x}}  \mathbb{P}\left( U_{\mathbf{x}} \geq snp^{\frac{1}{d}}\right).
\]
Given $\mathbf{x}$, choose any set $S_\mathbf{x}$ of edges of cardinality
\[
\left\lfloor \Cr{c: appendix_1} p^{\frac{1}{d} - 1}n \right\rfloor - \#\left\{e \in \cup_{i=0}^r (\ell x_i + [-\ell,\ell]^d)\right\}
\]
which does not intersect $\cup_{i=0}^r (\ell x_i + [-\ell,\ell]^d)$. Defining
\[
V_\mathbf{x} = U_\mathbf{x} + \sum_{f \in S_\mathbf{x}} X_f,
\]
we have $U_{\mathbf{x}} \leq V_{\mathbf{x}}$, $\mathbb{E}V_{\mathbf{x}} \leq \Cr{c: appendix_1} p^{1/d}n$, and therefore
\begin{equation}\label{eq: overall_new_N_n}
\mathbb{P}\left( N_n \geq snp^{\frac{1}{d}}\right) \leq \sum_{\mathbf{x}}  \mathbb{P}\left( V_{\mathbf{x}} - \mathbb{E}V_{\mathbf{x}} \geq \left( \frac{s}{\Cr{c: appendix_1}} - 1 \right) \Cr{c: appendix_1} np^{\frac{1}{d}}\right).
\end{equation}
To estimate this probability, we use the following lemma, which is \cite[Cor.~2.4]{J02}. For its statement, let $k' = \#\{ f \in \mathcal{E}^d : d(e,f) \leq k\}$ (for any $e \in \mathcal{E}^d$).

\begin{lem}
Let $e_1, \dots, e_m \in \mathcal{E}^d$ be distinct and define $\phi(x) = (1+x)\log (1+x)-x$ for $x\geq 0$. For all $t > 0$, we have
\begin{equation}\label{eq: concentration_1}
\mathbb{P}\left( \sum_{i=1}^m X_{e_i} - mp > t\right) \leq \exp\left( - \frac{mp}{(k'+1)(1-p)} \phi\left( \frac{4t}{5mp}\right)\right).
\end{equation}
\end{lem}

Applying \eqref{eq: concentration_1} with $m = \lfloor \Cr{c: appendix_1}p^{(1/d)-1}n\rfloor$ and using the fact that $\phi(x)$ is nondecreasing for $x \geq 0$, we obtain
\begin{align*}
\mathbb{P}\left( V_{\mathbf{x}} - \mathbb{E}V_{\mathbf{x}} \geq \left( \frac{s}{\Cr{c: appendix_1}} - 1 \right) \Cr{c: appendix_1} np^{\frac{1}{d}}\right) &\leq \exp\left( - \frac{\left\lfloor \Cr{c: appendix_1} p^{\frac{1}{d} - 1}n \right\rfloor p}{(k'+1)(1-p)}  \phi\left( \frac{4 \left(\frac{s}{\Cr{c: appendix_1}} - 1 \right) \Cr{c: appendix_1} np^{\frac{1}{d}}}{5\left\lfloor \Cr{c: appendix_1} p^{\frac{1}{d} - 1}n \right\rfloor p} \right)\right) \\
&\leq \exp\left(  - \frac{\Cr{c: appendix_1} p^{\frac{1}{d}}n}{2(k'+1)} \phi\left( \frac{4}{5} \left(\frac{s}{\Cr{c: appendix_1}} - 1\right) \right)\right).
\end{align*}
For some constant $\Cl[lgc]{c: appendix_2}>0$, we have $\phi((4/5)((s/\Cr{c: appendix_1})-1)) \geq s/2$ whenever $s \geq \Cr{c: appendix_2}$, so we can return to \eqref{eq: overall_new_N_n} to obtain
\begin{align*}
\mathbb{P}\left( N_n \geq snp^{\frac{1}{d}}\right) \leq \sum_{\mathbf{x}} \exp\left(  - \frac{\Cr{c: appendix_1} p^{\frac{1}{d}}n}{4(k'+1)} s\right) &\leq 3^{5dnp^{\frac{1}{d}}} \exp\left(  - \frac{\Cr{c: appendix_1} p^{\frac{1}{d}}n}{4(k'+1)} s\right) \\
&= \exp\left( - \frac{p^{\frac{1}{d}}n}{k'+1} \left( \frac{\Cr{c: appendix_1}s}{4} - 5d(k'+1)\log 3\right)\right).
\end{align*}
Here, we have used that the number of choices of $\mathbf{x}$ is at most $3^{dr} \leq 3^{5dnp^{1/d}}$. For $s \geq \Cl[lgc]{c: appendix_3}(k'+1)$, where $\Cr{c: appendix_3}>0$ is a constant, we have the bound
\begin{equation}\label{eq: my_new_equation_2024}
\mathbb{P}\left( N_n \geq snp^{\frac{1}{d}}\right) \leq \exp\left( - \frac{\Cr{c: appendix_1}p^{\frac{1}{d}}n}{8(k'+1)} s\right) \leq \exp\left( - \frac{\Cr{c: appendix_1}}{{8(k'+1)}} s\right).
\end{equation}
Last, we integrate this inequality to produce, for some $\Cr{c: appendix_4}>0$,
\[
\mathbb{E}\left( \frac{N_n}{np^{\frac{1}{d}}}\right) \leq \Cr{c: appendix_3}(k'+1) + \int_{\Cr{c: appendix_3}(k'+1)}^\infty \mathbb{P}(N_n \geq snp^{\frac{1}{d}})~\text{d}s \leq \Cl[lgc]{c: appendix_4}k^d.
\]
\end{proof}

The proof of Prop.~\ref{prop: lattice_animals} implies estimates for the $\ell$-th moment of $N_n$ with $\ell > 1$. A special case we use earlier in the paper is the following. There exists $\Cl[lgc]{c: 2024_constant_1}>0$ such that for all $p \in (0,1)$, all $n \geq 1$, and all families $(X_e)_{e \in \mathcal{E}^d}$ of i.i.d.~Bernoulli$(p)$ random variables,
\begin{equation}\label{eq: lattice_animals_second_moment}
\mathbb{E}N_n^2 \leq \left( \Cr{c: 2024_constant_1} p^{\frac{1}{d}} n\right)^2.
\end{equation}
\begin{proof}[Proof of~\eqref{eq: lattice_animals_second_moment}]
As in the proof of Prop.~\ref{prop: lattice_animals}, if $n \leq p^{-1/d}$, then using $\#\{e : e \in [-n,n]^d\} \leq \Cl[lgc]{c: 2024_constant_2}n^d$ for $\Cr{c: 2024_constant_2} = 3^{d+1}d$, we get
\begin{align*}
\mathbb{E}\left( \frac{N_n}{p^{\frac{1}{d}}n} \right)^2 \leq \frac{1}{\left( p^{\frac{1}{d}}n\right)^2} \mathbb{E}\left(\sum_{e\in [-n,n]^d}X_e\right)^2 &= \frac{1}{\left( p^{\frac{1}{d}}n\right)^2} \left( \Cr{c: 2024_constant_2}n^d p(1-p) + (\Cr{c: 2024_constant_2}n^d p)^2\right) \\
&\leq \Cr{c: 2024_constant_2} \left( n p^{\frac{1}{d}}\right)^{d-2} + \Cr{c: 2024_constant_2}^2 \left( np^{\frac{1}{d}} \right)^{2d-2} \\
&\leq \Cr{c: 2024_constant_2} + \Cr{c: 2024_constant_2}^2.
\end{align*}

Instead if $n > p^{-1/d}$, then \eqref{eq: my_new_equation_2024} gives for $s \geq 2\Cr{c: appendix_3}$ (with $k'=1$)
\[
\mathbb{P}\left( N_n \geq snp^{\frac{1}{d}}\right)  \leq \exp\left( - \frac{\Cr{c: appendix_1}}{16} s\right).
\]
Upon integrating, we obtain $\mathbb{E}(N_n/np^{1/d})^2 \leq \Cl[lgc]{c: 2024_constant_3}$.
\end{proof}

\bigskip
\noindent
{\bf Acknowledgements.} Part of this work was done at the 2019 AMS Mathematics Research Communities workshop ``Stochastic Spatial Models'' which was funded by NSF grant DMS-1916439.

\end{document}